\documentclass[a4paper, 12pt, reqno]{amsart}

\usepackage{amscd}
\usepackage{amsmath}
\usepackage{amssymb}
\usepackage{amsthm}
\usepackage{bm}
\usepackage{color}
\usepackage{colortbl}
\usepackage{comment}
\usepackage{graphics}
\usepackage{Here}
\usepackage[nosort]{cite}
\usepackage{tikz}
\usepackage{tikz-cd}
\usepackage{latexsym}
\makeatletter
\@addtoreset{equation}{section}

\makeatother
\usetikzlibrary{matrix}
\title{Calabi--Yau 3-folds from projective joins of del Pezzo manifolds}
\date{}
\author[D.~Inoue]{Daisuke Inoue}
\address{
Graduate School of Mathematical Sciences,
The University of Tokyo,
3-8-1 Komaba,
Meguro-ku,
Tokyo,
153-8914,
Japan.}
\email{ino@ms.u-tokyo.ac.jp}

\def\E{\mathcal{E}}
\def\F{\mathcal{F}}

\def\I{\mathrm{I}}
\def\J{\mathcal{J}}
\def\K{\mathcal{K}}

\def\cO{\mathcal{O}}

\def\PP{\mathbb{P}}
\def\CC{\mathbb{C}}
\def\RR{\mathbb{R}}
\def\QQ{\mathbb{Q}}
\def\ZZ{\mathbb{Z}}

\def\pl{\oplus}
\def\qed{\hfill $\Box$}
\def\rar{\rightarrow}

\def\rank{\mathrm{rank}}
\def\det{\mathrm{det}}

\def\Join{\mathrm{Join}}

\def\Bl{\mathrm{Bl}}
\def\Gr{G(2, V_5)}
\def\pt{\mathrm{pt}}

\def\Perf{\mathrm{Perf}}

\newtheorem{Thm}{Theorem}[section]
\newtheorem{Lem}[Thm]{Lemma}
\newtheorem{Prop}[Thm]{Proposition}
\theoremstyle{definition}

\newtheorem{Conj}[Thm]{Conjecture}

\newtheorem{Fact}[Thm]{Fact}

\newtheorem{Rem}[Thm]{Remark}

\newtheorem*{acknowledgements}{Acknowledgements}

\newtheorem*{Pf}{Proof}
\newtheorem*{PfMain}{Proof of Theorem \ref{X_Y_birat}}
\renewcommand{\subsubsection}{%
\arabic{subsection}}

\begin{document}
\begin{abstract}
In this paper, we will construct new examples of derived equivalent Calabi--Yau $3$-folds with Picard number greater than one. 
We also study their mirror Calabi--Yau manifolds and find that they are given by Schoen's fiber products of suitable rational elliptic modular surfaces. 
\end{abstract}

\maketitle
\section{Introduction}\label{sec1}
The derived equivalence between Grassmannian and Pfaffian Calabi--Yau $3$-folds \cite{BCFKvS, Rod} is one of the interesting phenomena discovered in the study of mirror symmetry of Calabi--Yau $3$-folds. 
These Calabi--Yau $3$-folds share the same mirror family where the derived equivalence is indicated by the appearance of the corresponding (maximal) degenerations of the family. 
In this paper, we will construct some more examples of derived equivalent Calabi--Yau $3$-folds with Picard number greater than one. 
Our construction is an application of homological projective duality (HPD) of projective joins proved in a recent paper by Kuznetsov--Perry \cite{KP}. 
\vspace{3mm}

For a Calabi--Yau $3$-fold $X$, let us recall the so-called Fourier--Mukai partners of $X$ defined by 
\begin{align*}
\mathrm{FM}(X) = \{Y \mid D^b(X) \cong D^b(Y) \} / \sim_{\text{bir}}
\end{align*}
where $\sim_{\text{bir}}$ is a birational equivalence. 
The numbers of equivalent classes $|\mathrm{FM}(X)|$ are of considerable interest from both birational geometry as well as mirror symmetry. 
The derived equivalence between Grassmannian and Pfaffian Calabi--Yau $3$-folds was the first non-trivial example of $|\mathrm{FM}(X)| > 1$ with Picard number $\rho(X) = 1$ \cite{BC, Kuz_Gr_line}. 
When $\rho(X) = 1$, some other interesting examples of $|\mathrm{FM}(X)| > 1$ have been discovered; 
for example, the derived equivalence related to the Reye congruence and double quintic symmetroid Calabi--Yau $3$-folds by \cite{HT_Reye}, intersections of two Grassmannians by \cite{KP}, $G_2$-Grassmann Calabi--Yau $3$-folds by \cite{IMOU, KuzIMOU}, and Kapustka--Rampazzo's Calabi--Yau $3$-folds by \cite{KR}. 

When $\rho(X) > 1$, several possibilities can happen. 
For example, $X$ may have birational models which are not deformation equivalent to $X$. 
Also, $X$ may have a fibration $\pi : X \rar S$ whose general fiber is a lower dimensional Calabi--Yau manifold. 
For such a fibration, Bridgeland--Maciocia \cite{BM} considered a relative moduli space of stable sheaves $\mathcal{M}^s(X/S)$ and proved derived equivalence between $X$ and some component of $\mathcal{M}^s(X/S)$. 
Based on the result in \cite{BM}, C\u ald\u araru \cite{Cal} gave a sufficient condition for $X$ with elliptic fibration to be $|\mathrm{FM}(X)| > 1$. %
Similarly, Schnell \cite{Schn} has given an example $X$ with $|\mathrm{FM}(X)| > 1$ by studying a certain Calabi--Yau $3$-fold with abelian surface fibration due to Gross--Popescu \cite{GP}. 
\vspace{3mm}

In this paper, we will study Calabi--Yau $3$-folds with $|\mathrm{FM}(X)| > 1$ and $\rho(X) > 1$ from the viewpoint of projective join and its categorical version \cite{KP}. 
\vspace{3mm}

We first construct Calabi--Yau $3$-folds as linear sections of projective joins of two projective varieties. 
Here we recall the following well-known fact (cf.\ \cite{GP}): 

\begin{Fact}\label{join_ellips}
Let $E_1 \hookrightarrow \PP(V_1)$ and $E_2 \hookrightarrow \PP(V_2)$ be projectively normal elliptic curves where $V_i$ are $n_i$-dimensional vector spaces for some $n_i$ ($i=1,2$).  
Let $\Join(E_1, E_2)$ be the projective join of $E_1, E_2$ (see (\ref{defjoin}) for the definition). 
Then $\Join(E_1, E_2)$ is normal and its dualizing sheaf is trivial. %
\end{Fact}

Based on this fact, we will construct our Calabi--Yau $3$-folds by smoothings of $\Join(E_1, E_2)$, where $E_1$ is an elliptic curve of degree $5$ in $\Gr$ and $E_2$ is one of three possible choices of elliptic curves given as linear sections of suitable del Pezzo manifolds (see Section \ref{The_case_del_Pezzo}). %
Combining this construction with the homological projective duality for categorical joins in \cite{KP}, we will naturally come to a pair of Calabi--Yau $3$-folds which are derived equivalent. 
By studying birational geometry of these Calabi--Yau $3$-folds, we obtain the following result: 

\begin{Thm}\label{main_theorem}
For each $\Join(E_1, E_2)$ as above, there exists a pair of Calabi--Yau $3$-folds which are derived equivalent, but not birationally equivalent. Both Calabi--Yau $3$-folds are realised as complete intersections in some projective bundles by relative hyperplane classes. 
\end{Thm}
\vspace{3mm}

We will also construct Calabi--Yau $3$-folds which are mirror dual to those constructed in Theorem \ref{main_theorem}. 

We remark that the standard methods to have mirror dual do not apply to Calabi--Yau $3$-folds constructed in Theorem \ref{main_theorem} since they are not complete intersections in Gorenstein Fano toric varieties. 
However they are given as complete intersections by linear sections of some projective bundles (Section \ref{The_case_del_Pezzo}). 
We note that the so-called $I$-functions of Gromov--Witten theory are known for such Calabi--Yau manifolds \cite{Bro}. 
Since it is known that these $I$-functions coincide with the period integrals of the expected mirror families of Calabi--Yau $3$-folds, we can try and error to find mirror duals of Calabi--Yau manifolds in Theorem \ref{main_theorem}. 
In fact, in the talk \cite{Gal}, Galkin remarked that the $I$-functions of Calabi--Yau $3$-folds in projective join of del Pezzo manifolds are given by {\it Hadamard products} of $I$-functions of elliptic curves which are linear sections of del Pezzo manifolds. 
We note that Hadamard products are used to construct interesting $4$-th order Fuchsian differential equations (e.g.\ \cite{AZ, vS}). 
Also, Hadamard products of period integrals of families of elliptic curves are interpreted as period integrals of certain families of Schoen's Calabi--Yau $3$-folds (see \cite{SvS}). 
Combining these, we will propose in general the following: 

\begin{Conj}
Mirror Calabi--Yau $3$-folds of linear section Calabi--Yau $3$-folds in projective join of del Pezzo manifolds are given by Schoen's Calabi--Yau $3$-folds with certain choices of rational elliptic surfaces with section. 
\end{Conj}

In this paper, we will provide supporting evidences of this conjecture by constructing two parameter families of Schoen's Calabi--Yau $3$-folds as the mirror manifolds of Calabi--Yau $3$-folds in Theorem \ref{main_theorem}. 
In our construction, it turns out that these family of Calabi--Yau $3$-folds have several large complex structure limit points. 
These large complex structure limit points can be interpreted which as non-trivial Fourier--Mukai partners, i.e., $|\mathrm{FM}(X)| > 1$, as well as birational models of $X$. 

Recently, Hosono and Takagi studied a relation between mirror symmetry and birational geometry \cite{HT}. 
Our examples of Calabi--Yau $3$-folds should also be interesting from their point of view. 
\vspace{3mm}

The organization of this paper is as follows. 
In Section \ref{sec2}, we will summarize projective joins. 
We will also introduce examples of pairs of varieties which are homological projective dual to each other. %
In Section \ref{sec3}, we will show the main result ($=$Theorem \ref{main_theorem}) in this paper. 
We will construct several pairs of Calabi--Yau $3$-folds from mutually orthogonal linear sections of suitable projective bundles. 
We will show that each pair of Calabi--Yau $3$-folds is not birational. 
In Section \ref{sec4}, we will construct Calabi--Yau $3$-folds which are topological mirror of Calabi--Yau $3$-folds constructed in Section \ref{sec3}. 
In Section \ref{sec5}, we will construct two dimensional families of mirror Calabi--Yau $3$-folds. 
We will see that this family of Calabi--Yau $3$-folds have several large complex structure limit points. %

\begin{acknowledgements}
The author is deeply grateful to Shinobu Hosono for many suggestions and helpful comments. 
The author is also grateful to Hiromichi Takagi for valuable comments to improve this paper. 
The author would like to thank Atsushi Kanazawa for encouragement during the postdoctoral period. 
\end{acknowledgements}

\section{Joins and homological projective duality for projective bundles}\label{sec2}
In this section, following \cite{KP}, we introduce projective joins and resolved joins of two projective varieties. %
We also introduce pairs of projective varieties which are homological projective dual to each other. 

\subsection{Projective joins}
Let $M_1 \subset \PP(V_1)$ and $M_2 \subset \PP(V_2)$ be smooth projective varieties. 
The projective join of $M_1$ and $M_2$ is defined by 
\begin{equation}\label{defjoin}
\Join(M_1, M_2) := \bigcup_{x_1 \in M_1, x_2 \in M_2}\langle [x_1,0], [0,x_2] \rangle
\end{equation}
where $\langle [x_1,0], [0,x_2] \rangle$ is the line spanned by the points $[x_1,0]$ and $[0,x_2]$ in $\PP(V_1 \pl V_2)$. %
By the natural embedding, we can regard $M_i$ as a subvariety of $\Join(M_1, M_2)$. 
In general, $\Join(M_1, M_2)$ is singular along the disjoint union $M_1 \sqcup M_2$. 
There exists a resolution of $\Join(M_1, M_2)$, so-called {\it resolved join}, which is given by 
\begin{equation}\label{resjoin}
\PP_{M_1 \times M_2}(\mathcal{O}(-H_1) \pl \mathcal{O}(-H_2))
\end{equation}
where $H_i$ is the hyperplane class on $M_i$ for $i=1,2$. 
For a projective variety $M$, we denote by $M^{\ast}$ the projective dual variety of $M$. 
The following formula is well-known: 
\begin{equation}\label{clprojdual}
\Join(M_1, M_2)^{\ast} = \Join(M_1^{\ast}, M_2^{\ast}). 
\end{equation}
In a recent paper, Kuznetsov and Perry defined a certain category associated to the projective join (\ref{defjoin}), and introduced categorical join as a categorical extension of the classical duality (\ref{clprojdual}) in the setting of homological projective duality. 

\subsection{Examples of homological projective dual pairs}
By taking a linear section of projective joins (\ref{defjoin}), we will construct Calabi--Yau $3$-folds which are related to $\Join(E_1, E_2)$ for elliptic curves $E_i$ in $\PP(V_i)$. 
To apply the homological projective duality for categorical joins to our Calabi--Yau $3$-folds, we need to find homological projective dual varieties $X$ and $X^{\natural}$ together with morphisms $X \rar \PP(V)$ and $X^{\natural} \rar \PP(V^{\vee})$ such that mutually orthogonal linear sections 
\begin{equation*}\label{HPDellip}
X \times_{\PP(V)} \PP(W) \hspace{3mm} \text{and} \hspace{3mm} X^{\natural} \times_{\PP(V^{\vee})} \PP(W^{\perp})
\end{equation*}
are elliptic curves where $W \subset V$ is a general linear subspace with an appropriate dimension, $V^{\vee}$ is the dual projective space of $V$, and $W^{\perp}$ is the orthogonal linear subspace to $W$. %
\vspace{3mm}

The most important example of $X$ and $X^{\natural}$ comes from $G(2, V_5)$, the Grassmannian of $2$-dimensional subspaces of $V_5 \cong \CC^5$. 
Let $G(2, V_5) \rar \PP(\wedge^2 V_5)$ be the Pl\"{u}cker embedding of $G(2, V_5)$. 
It is known \cite{KuzHyp} that the homological projective dual of the $G(2, V_5)$ with $G(2, V_5) \rar \PP(\wedge^2 V_5)$ is the dual Grassmannian $G(2, V_5^{\vee})$ with $G(2, V_5^{\vee}) \rar \PP(\wedge^2 V_5^{\vee})$. 
Let $W \subset \wedge^2 V_5$ be a general $5$-dimensional linear subspace and $W^{\perp} \subset \wedge^2 V_5^{\vee}$ be the orthogonal linear subspace to $W$. 
Then, both orthogonal linear sections 
\begin{equation*}
G(2, V_5) \times_{\PP(\wedge^2 V_5)} \PP(W) \hspace{3mm} \text{and} \hspace{3mm} G(2, V_5^{\vee}) \times_{\PP(\wedge^2 V_5^{\vee})} \PP(W^{\perp})
\end{equation*}
are elliptic curves of degree $5$. 
\vspace{3mm}

Other examples of $X$ and $X^{\natural}$ may be given by the linear duality \cite{KuzHPD}, which is the homological projective duality valid for projective bundles. 
The following construction gives Calabi--Yau $n$-folds as linear sections of projective bundles over an $(n+1)$-dimensional Fano manifold $Z$. 
This is a special case of \cite[Section 8]{KuzHPD} when we take $\F = \cO_Z^{\oplus i}$ in the notation there. %

\makeatletter
\renewcommand{\theenumi}{(\@roman\c@enumi)}
\renewcommand{\labelenumi}{\theenumi}
\makeatother

\begin{Prop}\label{FanotoCYn}
Let $Z$ be an $(n+1)$-dimensional Fano manifold. 
Let $\E$ be a locally free sheaf on $Z$ of rank $r$. 
Let $\PP(\E)$ be the projective bundle over $Z$ associated to $\E$. 
We denote by $\pi$ the projection of $\PP(\E)$ to $Z$. 
We denote by $L$ the relative hyperplane class of $\PP(\E)$. 
Assume the following conditions: 
\begin{enumerate}
\item \label{ass(i)} the dual $\E^{\vee}$ is globally generated. 
\item \label{ass(ii)} $\dim \varphi_{|L|}(\PP(\E)) \geq r$ where $\varphi_{|L|}$ is the morphism defined by $|L|$. 
\item \label{ass(iii)} $c_1(\E) = K_Z$. 
\end{enumerate}
Let $X$ be a general complete intersection of $r$ relative hyperplanes in $\PP(\E)$. 
Then $X$ is a Calabi--Yau $n$-fold, i.e., $X$ is an $n$-dimensional smooth projective variety satisfying $K_X \sim 0$ and $h^i(\cO_X) = 0$ if $0 < i < n$ and $h^i(\cO_X) = 1$ if $i=0, n$. 
\end{Prop}
\begin{Pf}
Since the linear system $|L|$ is base point free by the assumption \ref{ass(i)}, by the Bertini theorem, a general complete intersection of $r$ divisors of $|L|$ is smooth. 
By the assumption \ref{ass(ii)}, we have 
\begin{align*}
\dim X = \dim \PP(\E) - r = n. 
\end{align*}
The canonical divisor of $\PP(\E)$ is 
\begin{align*}
K_{\PP(\E)} = -r L - \det \E + K_Z = -r L. %
\end{align*}
Hence, we see that $K_X \sim 0$ by the adjunction formula. 
To compute of dimensions $h^i(\cO_X) \; (0 \leq i \leq n)$, we use the Koszul resolution 
\begin{align*}
0 \rar \cO_{\PP(\E)}(-r L) \rar \cdots \rar \cO_{\PP(\E)}(-2L)^{\oplus \binom{r}{2}} \rar \cO_{\PP(\E)}(-L)^{\oplus r} \rar \cO_{\PP(\E)} \rar \cO_X \rar 0
\end{align*}
and the spectral sequence 
\begin{align}\label{spectral_Kos}
E_1^{-p, q} = H^q(\PP(\E), \cO_{\PP(\E)}(-p L)^{\oplus \binom{r}{p}}) \Rightarrow H^{q-p}(X, \cO_X). 
\end{align}
The $E_1$-term of (\ref{spectral_Kos}) can be computed by the spectral sequence 
\begin{align}\label{spectral_proj}
E_1^{p,q} = H^q(Z, R^p \pi_{\ast} \cO_{\PP(\E)}(-i L)) \Rightarrow H^{p+q}(\PP(\E), \cO_{\PP(\E)}(-i L))
\end{align}
for $i=0, \dots, r$. 
We can determine the higher direct images $R^p \pi_{\ast} \cO_{\PP(\E)}(-i L)$ by the Bott--Borel--Weil theorem. 
It turns that the following are the only cases where $R^p \pi_{\ast} \cO_{\PP(\E)}(-i L)$ are non-trivial: 
\begin{align*}
R^0 \pi_{\ast} \cO_{\PP(\E)} = \cO_Z, \hspace{3mm} R^{r-1} \pi_{\ast} \cO_{\PP(\E)}(-r L) = \cO_{Z}(K_Z). 
\end{align*}
Combining the spectral sequences (\ref{spectral_Kos}), (\ref{spectral_proj}) and the Kodaira vanishing on $Z$, we obtain the claimed results for $h^i(\cO_X)$. 
\qed
\end{Pf}

For a locally free sheaf $\E$ as in Proposition \ref{FanotoCYn}, we consider the orthogonal vector bundle of $\E$ which is defined by the kernel of the natural surjective homomorphism 
\begin{align}\label{exactseq_Eperp}
0 \rar \E^{\perp} \rar H^0(Z, \E^{\vee}) \otimes \cO_{Z} \rar \E^{\vee} \rar 0. 
\end{align}
From the exact sequence (\ref{exactseq_Eperp}), we see that $(\E^{\perp})^{\vee}$ is globally generated and $c_1(\E^{\perp}) = K_Z$. 
We denote by $r'$ the rank of $\E^{\perp}$. 
We denote by $L'$ the relative hyperplane class of $\PP(\E^{\perp})$. 
For simplicity, we now make the further assumptions that 
\begin{enumerate}
\setcounter{enumi}{3}
\item \label{ass(iv)} $\dim \varphi_{|L'|}(\PP(\E^{\perp})) \geq r'$. 
\item \label{ass(v)} the natural homomorphism $H^0(Z, \E^{\vee})^{\vee} \rar H^0(Z, (\E^{\perp})^{\vee})$ is an isomorphism. In other words, $H^1(Z, \E) = 0$. 
\end{enumerate}
Applying Proposition \ref{FanotoCYn} to $\PP(\E^{\perp})$, we have a Calabi--Yau $n$-fold as a general complete intersection of $r'$ relative hyperplanes in $\PP(\E^{\perp})$. 

Let us define $V = H^0(Z, \E^{\vee})^{\vee}$. 
The homological projective dual of the projective bundle $\PP(\E)$ with $\varphi_{|L|} : \PP(\E) \rar \PP(V)$ is the projective bundle $\PP(\E^{\perp})$ with $\varphi_{|L'|} : \PP(\E^{\perp}) \rar \PP(V^{\vee})$. 
This result is called linear HPD proved by Kuznetsov. 
As a consequence of HPD, we obtain the following result. 
\begin{Prop}
Let $\E$ be a locally free sheaf of rank $r$ on $Z$ satisfying the assumptions \ref{ass(i)}--\ref{ass(v)}. 
Let $W$ be a general codimension $r$ linear subspace of $V$. 
Let $W^{\perp}$ be an orthogonal linear subspace to $W$. 
Let $X = \PP(\E) \times_{\PP(V)} \PP(W)$ and $Y = \PP(\E^{\perp}) \times_{\PP(V^{\vee})} \PP(W^{\perp})$ be linear sections of projective bundles. 
Then $X$ and $Y$ are $n$-dimensional Calabi--Yau manifolds. 
Moreover, there exists an equivalence 
\begin{align*}
D^b(X) \cong D^b(Y)
\end{align*}
of the derived categories of these manifolds. 
\end{Prop}

\begin{Rem}
Let $D$ be the image of $X$ under the projection $\pi : \PP(\E) \rar Z$. 
This is equal to the image of $Y$ under the projection $\pi' : \PP(\E^{\perp}) \rar Z$, so we have the following diagram: 
$$
	\begin{tikzpicture}
	\draw[thick,->] (1.9,1.5) -- (0.2,1);
	\draw[thick,->] (1.9,0.5) -- (0.2,0);
	\draw[thick,->] (-1.9,1.5) -- (-0.2,1);
	\draw[thick,->] (-1.9,0.5) -- (-0.2,0);
	\node[above] at (2,1.5) {$\PP(\E^{\perp})$};
	\node[above] at (2,1) {$\cup$};
	\node[above] at (2,0.5) {$Y$};
	\node[below] at (0,1) {$Z$};
	\node[below] at (0,0.5) {$\cup$};
	\node[below] at (0,0) {$D$};
	\node[above] at (-2,1.5) {$\PP(\E)$};
	\node[above] at (-2,1) {$\cup$};
	\node[above] at (-2,0.5) {$X$}; 
	\end{tikzpicture} %
	$$ 
If $s_1, \dots, s_r \in H^0(Z, \E^{\vee})$ are defining sections of $X$, then the image $D$ is the anti-canonical hypersurface defined by $s_1 \wedge \cdots \wedge s_r \in H^0(Z, \wedge^r \E^{\vee}) = H^0(Z, \cO_Z(-K_Z))$. 
Therefore we have a birational map $X \dashrightarrow Y$. 
\end{Rem}

In the next section, we will use the construction of a Calabi--Yau $n$-fold in Proposition \ref{FanotoCYn} in the case that $S$ is a del Pezzo surface.

\section{Calabi--Yau $3$-folds from projective joins}\label{sec3}
In this section, we will construct three examples of Calabi--Yau $3$-folds, $X_i \; (i=1,2,3)$, as linear sections of resolved joins (\ref{resjoin}). 
For each of these Calabi--Yau $3$-folds, we obtain a Calabi--Yau $3$-fold $Y_i$ which is derived equivalent to $X_i$. 
By studying birational geometry of these Calabi--Yau $3$-folds, we will find that $X_i$ and $Y_i$ are not birational, hence non-trivial Fourier--Mukai partners to each other. 

\subsection{Joins of $G(2, V_5)$ and projective bundles over del Pezzo surfaces}
Let $M_1 = G(2, V_5)$ and $M_2$ be a projective bundle $\PP_S(\E)$ over a del Pezzo surface $S$ satisfying the assumptions \ref{ass(i)}--\ref{ass(v)}. 
We denote by $\Sigma_1$ the image of $M_1$ under the Pl\"{u}cker embedding. 
We also denote by $\Sigma_2$ the image of $M_2$ under the map corresponding to the relative hyperplane class. %

\begin{Prop}\label{CY3from_join}
Let $\PP_{M_1, M_2} = \PP_{M_1 \times M_2}(\cO(-H_1) \pl \cO(-H_2))$ be the resolved join of $M_1$ and $M_2$, where $H_1$ is the Schubert divisor class of $M_1$ and $H_2$ is the relative hyperplane class of $M_2$. 
Let $L$ be the relative hyperplane class of $\PP_{M_1, M_2}$. 
Then a general complete intersection of $r+5$ divisors in $|L|$ is a Calabi--Yau $3$-fold. 
\end{Prop}

\begin{Pf}
We will write $k=r+5$ here and hereafter in this section. 
Let us denote by $X$ a general complete intersection of $k$ divisors of $|L|$ in $\PP_{M_1, M_2}$. 

Since $H^0(\PP_{M_1, M_2}, \cO(L)) = H^0(M_1, \cO(H_1)) \pl H^0(M_2, \cO(H_2))$ and $H_i \; (i=1,2)$ is base point free by the assumption \ref{ass(i)}, the relative hyperplane class $L$ is also base point free. 
Then a general complete intersection $X$ is smooth by Bertini theorem. 

The image of $\PP_{M_1, M_2}$ under the map $\varphi_{|L|}$ is $\Join(\Sigma_1, \Sigma_2)$, which is the projective join of $\Sigma_1$ and $\Sigma_2$. 
By the assumption \ref{ass(ii)}, we have 
\begin{align*}
\dim \Join(\Sigma_1, \Sigma_2) = \dim \Sigma_1 + \dim \Sigma_2 + 1 \geq k. 
\end{align*}
Then we obtain $\dim X = \dim \PP_{M_1, M_2} - k = 3$. 

Similar to the proof of Proposition \ref{FanotoCYn}, we have 
\begin{align}\label{KX_resolved_join}
K_X = ((r + 3) L - 4H_1 - (r-1) H_2)|_X. 
\end{align}
Let $D_i$ be the divisor in $\PP_{M_1, M_2}$ corresponding to the injection $\cO(-H_i) \hookrightarrow \cO(-H_1) \pl \cO(-H_2)$ for $i=1,2$. 
The image of $D_i$ under the map $\varphi_{|L|}$ is $\Sigma_i \subset \Join(\Sigma_1, \Sigma_2)$. 
Since $\dim \Sigma_2 < k$, a general $X$ does not intersect with $D_2$. 
Similarly if $r \geq 2$, then a general $X$ does not intersect with $D_1$. 
These results imply that 
\begin{align*}
(L - H_i)|_X \sim 0 \hspace{3mm} (i=1,2). 
\end{align*}
By (\ref{KX_resolved_join}), we have 
\begin{align*}
K_X = 4(L - H_1)|_X + (r-1)(L - H_2)|_X \sim 0
\end{align*}
when the rank of $\E$ is greater than or equal to two. 
Similar result holds for the case $r=1$. 

To compute of dimensions $h^i(\cO_X) \; (0 \leq i \leq 3)$, we use the spectral sequence
\begin{align*}
E_1^{-p,q} = H^{q}(\PP_{M_1, M_2}, \cO_{\PP_{M_1, M_2}}(-pL))^{\oplus \binom{k}{p}} \Rightarrow H^{q-p}(X, \cO_X)
\end{align*}
associated with the Koszul resolution similar to the proof of Proposition \ref{FanotoCYn}. 
We note the following results which are obtained easily (e.g.\ Bott--Borel--Weil theorem): 
\begin{align}
&H^i(M_1, \cO(-j H_1)) = 
\begin{cases}
\CC &(\text{if} \; (i,j) = (0,0) \; \text{or} \; (6,5))\\
0 &(\text{if} \; 0 \leq j \leq 5 \; \text{and} \; (i,j) \neq (0,0) \; \text{or} \; (6,5)), 
\end{cases} \label{coh_M_1}\\
&H^i(M_2, \cO(-j H_2)) = 
\begin{cases}
\CC &(\text{if} \; (i,j) = (0,0) \; \text{or} \; (r+1,r)) \\
0 &(\text{if} \; 0 \leq j \leq r \; \text{and} \; (i,j) \neq (0,0) \; \text{or} \; (r+1,r)). 
\end{cases}
\end{align}
We also note that the higher direct images $R^p \pi_{\ast} \cO_{\PP_{M_1, M_2}}(-iL) \; (i=0, \dots, k)$ are nonzero only for $R^0 \pi_{\ast}\cO_{\PP_{M_1, M_2}} = \cO_{M_1 \times M_2}$ and 
\begin{align}\label{higher_direct_M_1_M_2}
R^1 \pi_{\ast} \cO_{\PP_{M_1, M_2}}(-i L) = \bigoplus_{a=0}^{i-2} \cO_{M_1 \times M_2}(-(a+1)H_1+(a+1-i)H_2)
\end{align}
for $2 \leq i \leq k$. 
As in the proof of Proposition \ref{FanotoCYn}, these results (\ref{coh_M_1})--(\ref{higher_direct_M_1_M_2}) give $h^i(\cO_X) = 1$ if $i=0,3$ and $h^i(\cO_X) = 0$ if $i=1,2$, i.e., $X$ is a Calabi--Yau $3$-fold. 
\qed
\end{Pf}

\begin{Rem}
We remark the relation between Fact \ref{join_ellips} and Proposition \ref{CY3from_join}. 
We denote by $V_{M_1} = H^0(M_1, \cO(H_1))^{\vee} = \wedge^2 V_5$ and $V_{M_2} = H^0(M_2, \cO(H_2))^{\vee}$. 
Let $E_1 = G(2, V_5) \times_{\PP(\wedge^2 V_5)} \PP(W_1)$ be a general linear section of $G(2, V_5)$ by $W_1$, where codimension of $W_1 \subset \wedge^2 V_5$ is five. 
Also let $E_2 = M_2 \times_{\PP(V_{M_2})} \PP(W_2)$ be a general linear section of $M_2$ by $W_2$, where codimension of $W_2 \subset V_{M_2}$ is $r$. 

Assume that the restriction of the map $M_2 \rar \Sigma_2$ to $E_2$ induces an isomorphism onto the image for a general choice of $W_2$. %
By Proposition \ref{CY3from_join}, for a general codimension $r + 5$ linear subspace $W \subset V_{M_1} \pl V_{M_2}$, the linear section $X = \PP_{M_1 \times M_2}(\cO(-H_1) \pl \cO(-H_2)) \times_{\PP(V_{M_1} \pl V_{M_2})} \PP(W)$ is a smooth Calabi--Yau $3$-fold. 
We define $\bar{X} = \varphi_{|L|}(X)$, which is a linear section of $\Join(\Sigma_1, \Sigma_2)$. 
If we consider a special choice $W = W_1 \pl W_2$ where $W_i \subset V_{M_i} \; (i=1,2)$ are given as above, then $\bar{X}$ is $\Join(E_1, E_2)$, which is nothing but a projective join of elliptic curves. 
\end{Rem}

\subsection{The cases where $\PP_S(\E)$ is a del Pezzo manifold}\label{The_case_del_Pezzo}
We will restrict our attention to the cases where a polarized manifold $(\PP_S(\E), L)$ is a del Pezzo manifold. 
From the classification of del Pezzo manifolds by Fujita and Iskovskikh (cf.\ \cite{Fuj, IP}), we consider the following three cases: 
\begin{align*}
(1)& \hspace{3mm} \PP^2 \times \PP^2 = \PP_{\PP^2}(\cO_{\PP^2}(-1)^{\pl 3}), \\
(2)& \hspace{3mm} \mathrm{Bl}_{\mathrm{pt}}\PP^3 = \PP_{\PP^2}(\cO_{\PP^2}(-2) \pl \cO_{\PP^2}(-1)), \\
(3)& \hspace{3mm} \PP^1 \times \PP^1 \times \PP^1 = \PP_{\PP^1 \times \PP^1}(\cO_{\PP^1 \times \PP^1}(-1, -1)^{\pl 2}). 
\end{align*}
We will denote by $N_i = \PP_{S_i}(\E_i) \; (i=1,2,3)$ the del Pezzo manifold of (1)--(3). 
Note that $\E_i$ satisfies the conditions \ref{ass(i)}--\ref{ass(v)}. %
We will denote by $r_i$ the rank of $\E_i$. 
We define $V_{N_i} = H^0(N_i, \cO(L_i))^{\vee}$ where $L_i$ is the relative hyperplane class on $N_i$. 
For each of these, we consider the following projective bundles associated to the orthogonal vector bundles %
\begin{align*}
(1)'& \hspace{3mm} \PP_{\PP^2}(\K_1^{\pl 3}), \\
(2)'& \hspace{3mm} \PP_{\PP^2}(\K_2 \pl \K_1), \\
(3)'& \hspace{3mm} \PP_{\PP^1 \times \PP^1}(\K_{1,1}^{\pl 2}), 
\end{align*}
by introducing locally free sheaves $\K_i \; (i=1,2)$, and $\K_{1,1}$ defined by the following exact sequences: 
\begin{align*}
&0 \rar \K_i \rar H^0(\PP^2, \cO_{\PP^2}(i)) \otimes \cO_{\PP^2} \rar \cO_{\PP^2}(i) \rar 0 \; (i=1,2), \\
&0 \rar \K_{1,1} \rar H^0(\PP^1 \times \PP^1, \cO_{\PP^1 \times \PP^1}(1,1)) \otimes \cO_{\PP^1 \times \PP^1} \rar \cO_{\PP^1 \times \PP^1}(1,1) \rar 0. 
\end{align*}
We will denote by $N'_i = \PP_{S_i}(\E_i^{\perp}) \; (i=1,2,3)$ the projective bundle of $(1)'$--$(3)'$. 
We will denote by $r_i'$ the rank of $\E_i^{\perp}$. 
\vspace{3mm}

Let $\mathcal{J}_i = \PP_{M_1 \times M_2}(\cO(-H_1) \pl \cO(-H_2))$ be the resolved join associated with $M_1 = G(2, V_5)$ and $M_2 = N_i$. 
Note that the maps $\varphi_{|L|} : \mathcal{J}_i \rar J_i$ are, respectively, the resolutions of the following projective joins: 
\begin{align*}
J_1 &= \Join(G(2, V_5), \PP^2 \times \PP^2), \\
J_2 &= \Join(G(2, V_5), \mathrm{Bl}_{\mathrm{pt}}\PP^3), \\
J_3 &= \Join(G(2, V_5), \PP^1 \times \PP^1 \times \PP^1) . 
\end{align*}

Let $X_i$ be a general complete intersection of $k_i:=r_i + 5$ divisors of $|L|$ in $\mathcal{J}_i$. 
In other words, we define $X_i$ by $\mathcal{J}_i \times_{\PP(\wedge^2 V_5 \pl V_{N_i})} \PP(W)$ for a general codimension $k_i$ subspace $W$ in $\wedge^2 V_5 \pl V_{N_i}$. 
By Proposition \ref{CY3from_join}, $X_i$ is a Calabi--Yau $3$-fold for a general $W$. 

\begin{Rem}
The restriction of $\varphi_{|L|}$ to $X_i$ induces an isomorphism $X_i$ to its image $\varphi_{|L|}(X_i)$ which is a linear section of the projective join of the corresponding del Pezzo manifolds. 
This construction of Calabi--Yau $3$-folds is due to S.\ Galkin \cite{Gal}. 
\end{Rem}

Let $\mathcal{J}_i' = \PP_{G(2, V_5^{\vee}) \times N_i'}(\cO(-H_1') \pl \cO(-H_2'))$ be the resolved join associated with $G(2, V_5^{\vee})$ and $N_i'$ where $H_1'$ is the Schubert divisor class of $G(2, V_5^{\vee})$ and $H_2'$ is the relative hyperplane class of $N_i'$. 
We denote by $L'$ the relative hyperplane class of $\J_i'$. 
We denote by $\Sigma'_i$ the image under the map $\varphi_{|H'_i|}$ for $i=1,2$. 
Since $H^0(\mathcal{J}_i', \cO(L'))$ is 
\begin{align*}
H^0(G(2, V_5^{\vee}), \cO(H_1')) \pl H^0(N_i', \cO(H_2')) = \wedge^2 V_5 \pl V_{N_i}
\end{align*}
by the assumption \ref{ass(v)}, a linear subspace $W$ in $\wedge^2 V_5 \pl V_{N_i}$ defines a linear section of $\mathcal{J}_i'$. 
In other words, we define $Y_i$ by 
\begin{align*}
Y_i = \mathcal{J}_i' \times_{\PP(\wedge^2 V_5^{\vee} \pl V_{N_i}^{\vee})} \PP(W^{\perp})
\end{align*}
where $W^{\perp}$ is the orthogonal subspace to $W$. 

\begin{Prop}\label{CY3_Hodge}
Let us fix a general linear subspace $W$ of codimension $k_i$ in $\wedge^2 V_5 \pl V_{N_i}$. 
Then the orthogonal linear section $Y_i$ is a Calabi--Yau $3$-fold. 
The Hodge numbers of these Calabi--Yau $3$-folds are given as in the following table: 
\begin{table}[H]
\centering
\begin{tabular}{rll}
\hline
 & $h^{1,1}$ & $h^{2,1}$ \\
\hline
$X_1$ & $2$ & $47$ \\
$X_2$ & $2$ & $47$ \\
$X_3$ & $3$ & $43$ \\
\hline
\end{tabular}
\hspace{5mm}
\begin{tabular}{rll}
\hline
 & $h^{1,1}$ & $h^{2,1}$ \\
\hline
$Y_1$ & $2$ & $47$ \\
$Y_2$ & $2$ & $47$ \\
$Y_3$ & $3$ & $43$ \\
\hline
\end{tabular}
\end{table}
\end{Prop}
\begin{Pf}
By definition, the equality $\dim V_{N_i} = r_i + r_i'$ holds. 
Then the codimension of $W^{\perp}$ in $\wedge^2 V_5^{\vee} \pl V_{N_i}^{\vee}$ is $r_i' + 5$. 
If we take a general codimension $r_i + 5$ linear subspace $W$ in $\wedge^2 V_{5} \pl V_{N_i}$, then $W^{\perp}$ is also a general codimension $r_i' + 5$ linear subspace in $\wedge^2 V_5^{\vee} \pl V_{N_i}^{\vee}$. 
By Proposition \ref{CY3from_join}, this implies that $Y_i$ is a Calabi--Yau $3$-fold. 

The calculations of the Hodge numbers of these Calabi--Yau $3$-folds are similar to the proof of Proposition \ref{CY3from_join}. 
For example, to compute the Hodge numbers of $X_1$, we use the conormal sequence
\begin{align*}
0 {\rar \cO(-L)^{\pl k_1}}|_{X_1} \rar \Omega_{\mathcal{J}_1}|_{X_1} \rar \Omega_{X_1} \rar 0. 
\end{align*}
The cohomologies of ${\cO(-L)^{\pl k_1}}|_{X_1}$ and $\Omega_{\mathcal{J}_1}|_{X_1}$ are calculated by the Koszul resolution and Bott--Borel--Weil theorem. 
Details are left to the readers. 
\qed
\end{Pf}

Applying the homological projective duality for categorical joins due to Kuznetsov and Perry, we obtain the following result. 
\begin{Prop}\label{KuzPer}
Assume that $X_i$ and $Y_i$ are given by mutually orthogonal linear sections in $\mathcal{J}_i$ and $\mathcal{J}_i'$. 
Then these Calabi--Yau $3$-folds are derived equivalent. 
\end{Prop}
\begin{Pf}
In the following, we only consider the case $i=1$. 
The other cases are proved similarly. 
\\ \;\;
Let $\mathcal{A}^1 = \Perf(G(2, V_5))$ and $\mathcal{A}^2 = \Perf(\PP_{\PP^2}(\cO(-1)^{\pl 3}))$ be the categories of perfect complexes. 
Let $(\mathcal{A}^1)^{\natural}$ and $(\mathcal{A}^2)^{\natural}$ be the HPD categories of these categories. 
There are equivalences
\begin{align*}
(\mathcal{A}^1)^{\natural} \cong \Perf(G(2,V_5^{\vee})), (\mathcal{A}^2)^{\natural} \cong \Perf(\PP_{\PP^2}(\K_1^{\pl 3}))
\end{align*}
as Lefschetz categories over $\PP(\wedge^2 V_5^{\vee})$ and $\PP(V_{N_1}^{\vee})$ by \cite[Corollary 8.3]{KuzHPD} and \cite{KuzHyp} respectively. 
Let $\J(\mathcal{A}^1, \mathcal{A}^2)$ be the categorical join of $\mathcal{A}^1$ and $\mathcal{A}^2$ defined by \cite[Definition 3.9]{KP}. 
By \cite[Theorem 4.1]{KP}, there exists an equivalence 
\begin{align*}
\J(\mathcal{A}^1, \mathcal{A}^2)^{\natural} \cong \J((\mathcal{A}^1)^{\natural}, (\mathcal{A}^2)^{\natural})
\end{align*}
of Lefschetz categories over $\PP(\wedge^2 V_5^{\vee} \pl V_{N_1}^{\vee})$. 
By the main theorem of HPD for joins (see \cite[Theorem 2.24]{KP}), we obtain an equivalence 
\begin{align}\label{equiv_lin_section}
\J(\mathcal{A}^1, \mathcal{A}^2)_{\PP(W)} \cong \J((\mathcal{A}^1)^{\natural}, (\mathcal{A}^2)^{\natural})_{\PP(W^{\perp})}
\end{align}
for a general codimension $k_1$ linear subspace $W$ of $\wedge^2 V_5 \pl V_{N_1}$. 

Similar to the proof of \cite[Lemma 6.13]{KP}, the left hand side of (\ref{equiv_lin_section}) is equivalent to $\Perf(X_1)$. 
Similarly, the right hand side of (\ref{equiv_lin_section}) is equivalent to $\Perf(Y_1)$. 
Therefore we have an equivalence 
\begin{align*}
\Perf(X_1) \cong \Perf(Y_1)
\end{align*}
for general $W$. 
By the argument of \cite[Remark 5.7]{KP}, this equivalence implies that 
\begin{align*}
D^b(X_1) \cong D^b(Y_1)
\end{align*}
as desired. 
\qed
\end{Pf}

\begin{Thm}\label{X_Y_birat}
Let $X_i$ and $Y_i$ be general orthogonal linear sections of $\mathcal{J}_i$ and $\mathcal{J}'_i$ for $1 \leq i \leq 3$. 
Then $X_i$ and $Y_i$ are not birationally equivalent. 
\end{Thm}

We prove this result by studying birational geometry of $X_i$ and $Y_i$. 
We give our proof only for the case $i=1$. 
The other cases are similar and left to the readers. %
\vspace{3mm}

Let $X_{1}$ be a linear section Calabi--Yau $3$-fold in $\J_1 = \PP_{\Gr \times \PP^2 \times \PP^2}(\cO(-H_1) \pl \cO(-H_2))$. 
We denote by $\pi_i$ the composition of the projection $\J_1$ to $\PP^2 \times \PP^2$ and the projections to the $i$-th factor for $i=1,2$. 

\begin{Lem}\label{X_1_ellip}
The restriction of $\pi_i$ to $X_1$ induces an elliptic fibration on $X_1$ for $i=1,2$. 
The divisors giving these elliptic fibrations generate the K\"{a}hler cone of $X_{1}$. 
In particular, the birational class of $X_1$ consists of $X_{1}$ itself. 
\end{Lem}

\begin{Pf}
The second assertion follows from the first assertion. 
In fact, $X_1$ has no flopping contraction because $X_{1}$ has Picard number two by Proposition \ref{CY3_Hodge}. 

Let us prove the first assertion. 
Let $D_i$ be the pull-back to $X_1$ of the hyperplane class of the $i$-th factor of $\PP^2 \times \PP^2$. 
The restriction of $\pi_1$ to $X_1$ induces a map $\pi_1 : X_1 \rar \PP^2$, which corresponds to the divisor $D_1$. 
This map gives an elliptic fibration on $X_1$. 
Indeed, a general fiber over a point $x \in \PP^2$ is 
\begin{align}\label{fiber_X_1}
\Join(G(2, V_5), \{ x \} \times \PP^2) \cap \PP(W)
\end{align}
where $W$ is the linear subspace of $\wedge^2 V_5 \pl V_{N_1}$ which defines $X_1$. 
Since $\Join(G(2, V_5), \{ x \} \times \PP^2)$ is a projective cone of $G(2, V_5)$ with its vertex $\PP^2$, the fiber (\ref{fiber_X_1}) is isomorphic to a linear section of $G(2, V_5)$ of codimension five, which is a degree $5$ elliptic curve. 
Similarly, the restriction of $\pi_2$ induces an elliptic fibration on $X_1$. 
\qed
\end{Pf}

Let us consider a linear section of $\mathcal{J}'_1 = \PP_{G(2, V_5^{\vee}) \times \PP_{\PP^2}(\K_1^{\pl 3})}(\cO(-H_1') \pl \cO(-H_2'))$. 
In the following, we will write the projective bundle $\PP_{\PP^2}(\K_1^{\pl 3})$ as $\PP_{\PP(V_3)}(\K_1 \otimes W_3^{\vee})$ or $\mathbb{P}_{\PP(W_3)}(V_3^{\vee} \otimes \K_1)$ where $V_3 \cong \CC^3$ and $W_3 \cong \CC^3$. 
Let us consider the following diagram: 
$$
	\begin{tikzpicture}
	\draw[thick,->] (2.9,1) -- (0.5,0) node[midway,below] {$\varphi_2$};
	\draw[thick,->] (-2.9,1) -- (-0.5,0) node[midway,below] {$\varphi_1$};
	\node[above] at (3,1) {$\mathbb{P}_{\PP(W_3)}(V_3^{\vee} \otimes \K_1)$};
	\node[below] at (0,0) {$\PP(V_3^{\vee} \otimes W_3^{\vee})$};
	\node[above] at (-3,1) {$\mathbb{P}_{\PP(V_3)}(\K_1 \otimes W_3^{\vee})$};
	\end{tikzpicture}
	$$
where $\varphi_i$ is the map defined by the relative hyperplane class for $i=1,2$. %
We use the following result in the appendix of \cite{HT}: 

\begin{Lem}\label{Spresol}
The maps $\mathbb{P}_{\PP(V_3)}(\K_1 \otimes W_3^{\vee}) \rar \PP(V_3^{\vee} \otimes W_3^{\vee})$ and $\mathbb{P}_{\PP(W_3)}(V_3^{\vee} \otimes \K_1) \rar \PP(V_3^{\vee} \otimes W_3^{\vee})$ are flopping contractions onto the common image
\begin{align*}
\Sigma'_2 = \{ M \in \PP(V_3^{\vee} \otimes W_3^{\vee}) \mid \rank \; M < 3\}
\end{align*}
where $\rank$ is a rank of a linear map $V_3 \rar W_3^{\vee}$ representing $M$. 
Moreover, the induced birational map $\mathbb{P}_{\PP(V_3)}(\K_1 \otimes W_3^{\vee}) \dashrightarrow \mathbb{P}_{\PP(W_3)}(V_3^{\vee} \otimes \K_1)$ is a flop. 
\end{Lem}

Note that if we identify $\PP(V_3^{\vee} \otimes W_3^{\vee})$ with the projectivization of the space of $3 \times 3$ matrices, then $\Sigma'_2$ is a cubic hypersurface defined by the determinant. 

Let $Y_{1}$ be the linear section Calabi--Yau $3$-fold in $\mathcal{J}'_1 = \PP_{G(2, V_5^{\vee}) \times \PP_{\PP(V_3)}(\K_1 \otimes W_3^{\vee})}(\cO(-H_1') \pl \cO(-H_2'))$ by a linear subspace $W' \subset \wedge^2 V_5^{\vee} \pl V_{N_1}^{\vee}$ of codimension $11$. 
We denote by $\bar{Y}_1$ the image of $Y_1$ under the map $\varphi_{|L'|}$ on $\mathcal{J}_1'$. 

\begin{Lem}\label{Y_1_small}
For general $W'$, the image $\bar{Y}_1$ has $30$ ODPs. %
Moreover, the restriction map $\varphi_{|L'|} : Y_1 \rar \bar{Y}_1$ is a small contraction. 
\end{Lem}
\begin{Pf}
The image of $\mathcal{J}'_1$ under the map $\varphi_{|L'|}$ is 
\begin{align*}
\Join(\Sigma'_1, \Sigma'_2) \subset \PP(\wedge^2 V_5^{\vee} \pl V_3^{\vee} \otimes W_3^{\vee}). %
\end{align*}
By using the map $\mathcal{J}'_1 \rar \Join(\Sigma'_1, \Sigma'_2)$, we see that the singular locus of $\Join(\Sigma'_1, \Sigma'_2)$ is equal to 
\begin{align*}
\Sigma'_1 \cup \Sigma'_2 \cup \Join(\Sigma'_1, D')
\end{align*}
where $D' \subset \Sigma'_2$ is a closed subvariety consisting of matrices with rank less than two. 
Note that $D'$ is isomorphic to $\PP^2 \times \PP^2$. 
Since $\bar{Y}_1$ is the linear section of $\Join(\Sigma'_1, \Sigma'_2)$ by a linear subspace of codimension $11$, $\bar{Y}_1$ does not intersect with the components $\Sigma'_1 \cup \Sigma'_2$. 
But $\bar{Y}_1$ must intersect with the component $\Join(\Sigma'_1, D')$ at finitely many points. 
These points are ODPs of $\bar{Y}_1$ because $\Sigma'_2$ has $A_1$-singularity along $D'$. %
Since $\deg \Join(\Sigma'_1, D') = 30$, $\bar{Y}_1$ has $30$ ODPs for general $W'$. 
The restriction map $\varphi_{|L'|} : Y_1 \rar \bar{Y}_1$ induces the isomorphism $\varphi_{|L'|}^{-1}(U) \cong U$ where $U = \bar{Y}_1 \setminus \mathrm{Sing}(\bar{Y}_1)$. 
Therefore the map $\varphi_{|L'|} : Y_1 \rar \bar{Y}_1$ is a crepant resolution, and the last claim follows from this fact. 
\qed
\end{Pf}

From the above diagram, we have the following diagram: 
$$
	\begin{tikzpicture}
	\draw[thick,->] (2.9,1) -- (0.5,0);
	\draw[thick,->] (-2.9,1) -- (-0.5,0);
	\node[above] at (4,1) {$\PP_{G(2, V_5^{\vee}) \times \PP_{\PP(W_3)}(V_3^{\vee} \otimes \K_1)}(\cO(-H_1') \pl \cO(-H_2'))$};
	\node[below] at (0,0) {$\Join(\Sigma'_1, \Sigma'_2)$};
	\node[above] at (-4,1) {$\PP_{G(2, V_5^{\vee}) \times \mathbb{P}_{\PP(V_3)}(\K_1 \otimes W_3^{\vee})}(\cO(-H_1') \pl \cO(-H_2'))$};
	\end{tikzpicture}
	$$ 
Combining Lemma \ref{Spresol} with Lemma \ref{Y_1_small}, we obtain the following result. 
\begin{Lem}
Let $\bar{Y}_1$ be a general linear section of $\Join(\Sigma'_1, \Sigma'_2)$. 
Let $Y_1$ and $Y_1'$ be the pull-backs of $\bar{Y}_1$ to $\PP_{G(2, V_5^{\vee}) \times \mathbb{P}_{\PP(V_3)}(\K_1 \otimes W_3^{\vee})}(\cO(-H_1') \pl \cO(-H_2'))$ and $\PP_{G(2, V_5^{\vee}) \times \PP_{\PP(W_3)}(V_3^{\vee} \otimes \K_1)}(\cO(-H_1') \pl \cO(-H_2'))$, respectively. 
Then both $Y_1$ and $Y_1'$ are Calabi--Yau $3$-folds and the maps $Y_1 \rar \bar{Y}_1$ and $Y_1' \rar \bar{Y}_1$ are small contractions. 
Moreover, the birational map $Y_1 \dashrightarrow Y_1'$ is a flop. %
\end{Lem}

Similar to the proof of Lemma \ref{X_1_ellip}, we have the following result. 
\begin{Lem}
There exist elliptic fibrations on $Y_1$ and $Y_1'$. 
\end{Lem}

The following diagram describes birational geometry of $Y_1$. 
$$
	\begin{tikzpicture}
	\draw[thick,->] (3.1,1) -- (5.5,0);
	\draw[thick,->] (2.9,2) -- (0.5,1);
	\draw[thick,->] (2.9,1) -- (0.5,0);
	\draw[thick,->] (-2.9,2) -- (-0.5,1);
	\draw[thick,->] (-2.9,1) -- (-0.5,0);
	\draw[thick,->] (-3.1,1) -- (-5.5,0);
	\node[above] at (4,2) {$\PP_{G(2, V_5^{\vee}) \times \PP_{\PP(W_3)}(V_3^{\vee} \otimes \K_1)}(\cO(-H_1') \pl \cO(-H_2'))$};
	\node[below] at (6,0) {$\PP^2$};
	\node[above] at (3,1.5) {$\cup$};
	\node[above] at (3,1) {$Y_1'$};
	\node[below] at (0,1) {$\Join(\Sigma'_1, \Sigma'_2)$};
	\node[below] at (0,0.5) {$\cup$};
	\node[below] at (0,0) {$\bar{Y}_1$};
	\node[above] at (-4,2) {$\PP_{G(2, V_5^{\vee}) \times \mathbb{P}_{\PP(V_3)}(\K_1 \otimes W_3^{\vee})}(\cO(-H_1') \pl \cO(-H_2'))$};
	\node[above] at (-3,1.5) {$\cup$};
	\node[above] at (-3,1) {$Y_1$};
	\node[below] at (-6,0) {$\PP^2$};
	\end{tikzpicture} %
	$$ 

\begin{PfMain}
Now, Theorem \ref{X_Y_birat} follows for $X_1$ and $Y_1$. 
Indeed, $X_1$ has no flopping type contraction, so $X_1$ and $Y_1$ cannot be birational. 
Other cases $i=2,3$ are proved by similar arguments. %
\qed
\end{PfMain}

\begin{Rem}
Both Calabi--Yau $3$-folds $X_1$ and $Y_1$ (or $Y_1'$) have elliptic fibrations over the same base $\PP^2$. 
It is interesting to see whether $Y_1$ (or $Y_1'$) is a component of the relative moduli space of stable sheaves on the elliptic fibrations $\pi_i : X_1 \rar \PP^2$ ($i=1,2$). 
\end{Rem}

In what follows, we summarize briefly the birational geometry of $X_2, Y_2$ and $X_3, Y_3$ in order. 
We first describe the birational geometry of $X_2$ and $Y_2$ in the following diagrams: 
$$
	
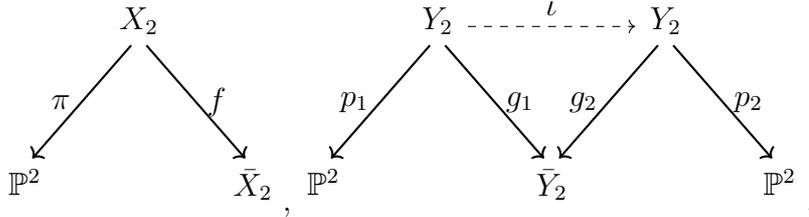
\begin{tikzpicture}
	\draw[thick,->] (0.1,1.5) -- (1.4,0) node[midway,right] {$f$};
	\draw[thick,->] (-0.1,1.5) -- (-1.4,0) node[midway,left] {$\pi$};
	\node[below] at (1.5,0) {$\bar{X}_2$};
	\node[above] at (0,1.5) {$X_2$};
	\node[below] at (-1.5,0) {$\PP^2$};
	\end{tikzpicture}, %
	\begin{tikzpicture}
	\draw[thick,->] (3.1,1.5) -- (4.4,0) node[midway,right] {$p_2$};
	\draw[thick,->] (2.9,1.5) -- (1.6,0) node[midway,left] {$g_2$};
	\draw[dashed,->] (0.4,1.75) -- (2.6,1.75) node[midway,above] {$\iota$};
	\draw[thick,->] (0.1,1.5) -- (1.4,0) node[midway,right] {$g_1$};
	\draw[thick,->] (-0.1,1.5) -- (-1.4,0) node[midway,left] {$p_1$};
	\node[below] at (4.5,0) {$\PP^2$};
	\node[above] at (3,1.5) {$Y_2$};
	\node[below] at (1.5,0) {$\bar{Y}_2$};
	\node[above] at (0,1.5) {$Y_2$};
	\node[below] at (-1.5,0) {$\PP^2$};
	\end{tikzpicture}. %
	$$
Here the map $\pi$ is an elliptic fibration on $X_2$ and $f$ is a divisorial contraction which contracts a degree $5$ del Pezzo surface to a point. 
The maps $p_i \; (i=1,2)$ are elliptic fibrations on $Y_2$ and $g_i \; (i = 1,2)$ are small contractions of $Y_2$. 
The rational map $\iota$ is a birational involution of $Y_2$ (see, for example, \cite[Proposition 6.1]{Ogu}). 

The following diagrams represent the sections of movable fans of $X_3$ and $Y_3$: 
$$
	
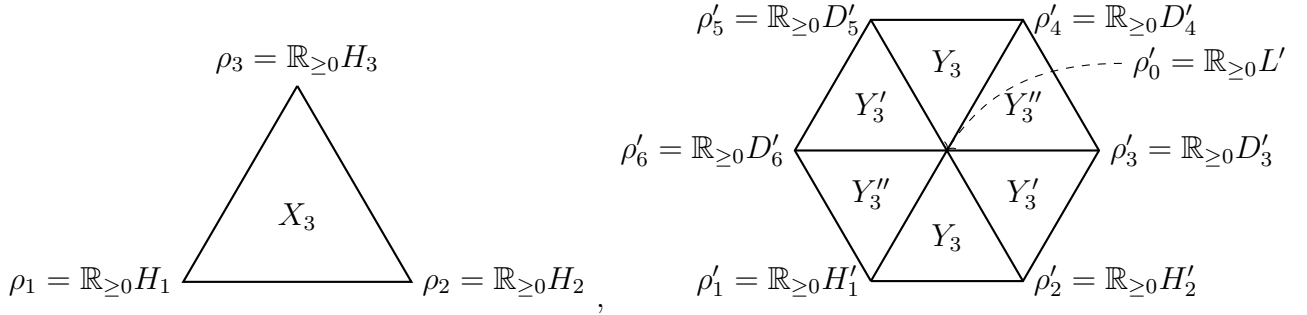
\begin{tikzpicture}
	\coordinate (A) at (0,{sqrt(3)*1.5-1.5});
	\coordinate (B) at (-1.5,-1.5);
	\coordinate (C) at (1.5,-1.5);
	\draw (0,{sqrt(3)*1.5/3-1.5}) node[circle] {$X_3$};
	\draw[thick] (A)-- (B)--(C)--(A);
	\node [above] at (A) {$\rho_3=\RR_{\geq 0}H_3$};
	\node [left] at (B) {$\rho_1=\RR_{\geq 0}H_1$};
	\node [right] at (C) {$\rho_2=\RR_{\geq 0}H_2$};
	\end{tikzpicture}, 
	\begin{tikzpicture}
	\coordinate (A) at (0,0);
	\coordinate (B) at (2,0);
	\coordinate (C) at (3,{sqrt(3)});
	\coordinate (D) at (2,{sqrt(3)*2});
	\coordinate (E) at (0,{sqrt(3)*2});
	\coordinate (F) at (-1,{sqrt(3)});
	\coordinate (G) at ({1},{sqrt(3)});
	\draw (1,{sqrt(3)/3}) node[circle] {$Y_3$};
	\draw (1,{5*sqrt(3)/3}) node[circle] {$Y_3$};
	\draw (2,{2*sqrt(3)/3}) node[circle] {$Y'_3$};
	\draw (2,{4*sqrt(3)/3}) node[circle] {$Y''_3$};
	\draw (0,{2*sqrt(3)/3}) node[circle] {$Y''_3$};
	\draw (0,{4*sqrt(3)/3}) node[circle] {$Y'_3$};
	\draw[thick] (A)-- (B)--(C)--(D)--(E)--(F)--(A);
	\draw[thick] (A)-- (G)--(B);
	\draw[thick] (C)-- (G)--(D);
	\draw[thick] (E)-- (G)--(F);
	\node [left] at (A) {$\rho'_1=\RR_{\geq 0} H_1'$};
	\node [right] at (B) {$\rho'_2=\RR_{\geq 0} H_2'$};
	\node [right] at (C) {$\rho'_3=\RR_{\geq 0} D_3'$};
	\node [right] at (D) {$\rho'_4=\RR_{\geq 0} D_4'$};
	\node [left] at (E) {$\rho'_5=\RR_{\geq 0} D_5'$};
	\node [left] at (F) {$\rho'_6=\RR_{\geq 0} D_6'$};
	\node [right] at (3.3,{5*sqrt(3)/3}) {$\rho'_0=\RR_{\geq 0}L'$};
	\draw[dashed,->] (3.3,{5*sqrt(3)/3}) to [out=180,in=58] (G);
	\end{tikzpicture}
	$$
where $Y'_3$ and $Y''_3$ are birational models of $Y_3$. %
The divisors $H_1, H_2, H_3$ (resp.\ $H_1', H_2', L'$) are the generator of the nef cone of $X_3$ (resp.\ $Y_3$). 
The divisors $D'_3, D'_4, D'_5, D'_6$ are $L'-H'_1+H'_2, 2L'-H_1', 2L'-H_2', L'+H_1'-H_2'$, respectively. %
We describe the morphisms associated to each cone of the movable fans: 
Two dimensional cones $\rho_i + \rho_j$ $(1 \leq i < j \leq 3)$ correspond to elliptic fibrations on $X_3$. 
One dimensional cones $\rho_i$ $(i=1,2,3)$ correspond to K3 fibrations on $X_3$. 
Similarly, two dimensional cones $\rho'_i + \rho'_{i+1}$ $(1 \leq i \leq 6)$ correspond to elliptic fibrations on some birational models $Y_3, Y'_3$, or $Y''_3$. 
Here we use the notation $\rho'_7 = \rho'_1$. 
One dimensional cones $\rho'_i$ $(1 \leq i \leq 6)$ correspond to K3 fibrations on some birational models $Y_3, Y'_3$, or $Y''_3$. 
Moreover, two dimensional cones $\rho'_0 + \rho'_i$ $(1 \leq i \leq 6)$ correspond to small contractions of birational models $Y_3, Y'_3$, or $Y''_3$. 
The one dimensional cone $\rho'_0$ corresponds to divisorial contractions of birational models $Y_3, Y'_3$, and $Y''_3$. 

\section{Fiber products of rational elliptic surfaces}\label{sec4}
In this section, we will construct Calabi--Yau $3$-folds which are mirror to $X_i \; (i=1,2,3)$ constructed in Section \ref{The_case_del_Pezzo}. 
In addition to these Calabi--Yau $3$-folds, we also consider a Calabi--Yau $3$-fold $X_0$ given as a linear section of the join $J_0 = \Join(\Gr, \Gr)$. 
This Calabi--Yau $3$-fold is known as intersections of two Grassmannians studied by \cite{GP, Kan, Kap_prim, Kap_tran}. 

The strategy of our construction is based on the following two observations. 
\vspace{3mm}

\begin{enumerate}
\item \label{ass4(i)} There is a relation between projective joins and Hadamard products of corresponding $I$-functions \cite{Gal}. %
\item \label{ass4(ii)} There is a geometric realization of a Hadamard product as a period integral of a fiber product of two families of elliptic curves. %
\vspace{3mm}
\end{enumerate}

These observations naturally lead us to a family of Calabi--Yau $3$-folds for each of our Calabi--Yau $3$-folds $X_i$ (or $Y_i$). 
To go into the details of mirror symmetry, we need to find suitable compactifications of the parameter spaces of the families. 
In this section, deferring this problem to the next section, we will restrict our attentions to constructing general members of such families and verifying the expected exchanges of the Hodge numbers. %

\subsection{Schoen's Calabi--Yau $3$-folds}\label{Schoen_CY}
Let us recall the work of Schoen \cite{Scho} to construct Calabi--Yau $3$-folds as crepant resolutions of fiber products of two rational elliptic surfaces with sections. 

Let $\pi_i : S_i \rar \PP^1$ $(i=1,2)$ be relatively minimal rational elliptic surfaces with sections. 
We denote by $D_i$ the image of singular fibers under the map $\pi_i$ for $i=1,2$. 
For simplicity, we assume that all singular fibers of $\pi_i$ $(i=1,2)$ over $D_1 \cap D_2$ are of type $\I_{b} \; (b \geq 0)$ in the notation of Kodaira. 
We set  
\begin{align*}
Y = S_1 \times_{\PP^1} S_2.  
\end{align*}
A point $(x_1,x_2) \in Y$ is singular if and only if $x_i$ is a singular point of the singular fiber $\pi_i^{-1}(t)$ for both $i=1, 2$ where $t := \pi_1(x_1) = \pi_2(x_2)$. 
By assumption, it is easy to check that these singular points are ordinary double points. 

Schoen proved that the dualizing sheaf of $Y$ is trivial. %
A crepant resolution of $Y$ is given by an iterated blow-up along smooth divisors on $Y$. 
Indeed, each node of $Y$ is lying on a smooth divisor of $Y$ which is a product of irreducible components of singular fibers $\pi_1^{-1}(t)$ and $\pi_2^{-1}(t)$ where $t \in D_1 \cap D_2$. 
We denote by $\tilde{Y}$ a crepant resolution of $Y$. 
Then $\tilde{Y}$ is a Calabi--Yau $3$-fold. 
Note that the condition $h^1(\cO_{\tilde{Y}}) = 0$ is proved by using the fact that $Y$ is a divisor in $S_1 \times S_2$. 

\subsection{Calabi--Yau $3$-folds from rational elliptic modular surfaces}\label{CY_from_modular}
We consider Schoen's Calabi--Yau $3$-folds in the case where $S_i$ are rational elliptic modular surfaces. 

From the table of elliptic surfaces with four singular fibers due to Herfurtner \cite[Table 3]{Her}, we take the elliptic surfaces $S(5), S(6), S(7)$ with the following types of singular fibers: 
\begin{table}[H]
\centering
\begin{tabular}{rllll}
\hline
$S(5)$ & $\I_5$ & $\I_5$ & $\I_1$ & $\I_1$ \\
$S(6)$ & $\I_6$ & $\I_3$ & $\I_2$ & $\I_1$ \\
$S(7)$ & $\I_7$ & $\I_2$ & $\I_1$ & $\I\I$ \\
\hline
\end{tabular}
\end{table}
These are relatively minimal rational elliptic surfaces with section. 
Note that $S(5)$ and $S(6)$ are Shioda's modular surfaces \cite{Shi} associated to the congruence subgroups of $SL(2, \ZZ)$ given by %
\begin{align*}
\Gamma_1(n) = 
\left\{
\left(
\begin{array}{cc}
a & b \\
c & d
\end{array}
\right) \in SL(2, \ZZ) \;
\middle| \; a, d \equiv 1, c \equiv 0 \; (\mathrm{mod} \; n)
\right\}
\end{align*}
with $n = 5,6$, respectively. %

By using an automorphism of $\PP^1$, we define four elliptic surfaces $T := T_{\bar{0}}, T_0, T_1, T_2, T_3$ with the following specified singular fibers at $0$ and $\infty$: 
\begin{itemize}
\item $T_{\bar{0}} \cong S(5)$, $D = \{0, \bar{x}_1, \bar{x}_2, \infty \}$, $\pi^{-1}(0)$ is of type $\I_5$, and $\pi^{-1}(\infty)$ is of type $\I_5$, 
\item $T_0 \cong S(5)$, $D = \{0, x_1, x_2, \infty \}$, $\pi^{-1}(0)$ is of type $\I_5$, and $\pi^{-1}(\infty)$ is of type $\I_5$, 
\item $T_1 \cong S(6)$, $D = \{0, y_1, y_2, \infty \}$, $\pi^{-1}(0)$ is of type $\I_6$, and $\pi^{-1}(\infty)$ is of type $\I_3$, 
\item $T_2 \cong S(7)$, $D = \{0, z_1, z_2, \infty \}$, $\pi^{-1}(0)$ is of type $\I_7$, and $\pi^{-1}(\infty)$ is of type $\I_2$, 
\item $T_3 \cong S(6)$, $D = \{0, w_1, w_2, \infty \}$, $\pi^{-1}(0)$ is of type $\I_6$, and $\pi^{-1}(\infty)$ is of type $\I_2$, %
\end{itemize}
where $\pi : T \rar \PP^1$ and $D$ is the image of singular fibers under the map $\pi$. 
We assume that singular fibers other than $0$ and $\infty$ are in general positions. 
Namely, all the points $\bar{x}_1, \bar{x}_2, x_1, x_2, y_1, y_2, z_1, z_2, w_1, w_2$ are mutually distinct. 
\vspace{3mm}

Based on the observations \ref{ass4(i)}, \ref{ass4(ii)} in the beginning of this section, we claim the following. 
\begin{Conj}\label{main_Conj}
Let $X_i^{\ast} = T_{\bar{0}} \times_{\PP^1} T_i$ $(i = 0,1,2,3)$ be the fiber products of relatively minimal rational elliptic surfaces with sections. 
Let $\tilde{X}_i^{\ast}$ $(i=0,1,2,3)$ be the Schoen's Calabi--Yau $3$-folds constructed in Section \ref{Schoen_CY}. 
Then $\tilde{X}_0^{\ast}$, $\tilde{X}_1^{\ast}$, $\tilde{X}_2^{\ast}$, and $\tilde{X}_3^{\ast}$, respectively, are mirror Calabi--Yau $3$-folds of the linear section Calabi--Yau $3$-folds of $\Join(G(2, V_5), G(2, V_5))$, $\Join(G(2, V_5), \PP^2 \times \PP^2)$, $\Join(G(2, V_5), \Bl_{\pt}\PP^3)$, and $\Join(G(2, V_5), \PP^1 \times \PP^1 \times \PP^1)$. 
\end{Conj}

\begin{Rem}
We can find these elliptic surfaces $T_{i}$ $(1 \leq i \leq 3)$ by taking special one parameter families of Batyrev--Borisov toric mirror constructions of complete intersection elliptic curves in $\PP^2 \times \PP^2, \Bl_{\pt}\PP^3$, and $\PP^1 \times \PP^1 \times \PP^1$, respectively. %
\end{Rem}

Below, we will provide convincing evidences for this conjecture. 

\subsection{Calculations of Hodge numbers}
We calculate Hodge numbers of $\tilde{X}_i^{\ast} \; (i=0,1,2,3)$ by the Schoen's results \cite{Scho}. 

\begin{Prop}
The Hodge numbers of $\tilde{X}_i^{\ast} \; (i=0,1,2,3)$ are given as in the following table: 
\begin{table}[H]
\centering
\begin{tabular}{rll}
\hline
 & $h^{1,1}$ & $h^{2,1}$ \\
\hline
$\tilde{X}_0^{\ast}$ & $51$ & $1$ \\
$\tilde{X}_1^{\ast}$ & $47$ & $2$ \\
\hline
\end{tabular}
\hspace{5mm}
\begin{tabular}{rll}
\hline
 & $h^{1,1}$ & $h^{2,1}$ \\
\hline
$\tilde{X}_2^{\ast}$ & $47$ & $2$ \\
$\tilde{X}_3^{\ast}$ & $43$ & $3$ \\
\hline
\end{tabular}
\end{table}
\noindent In particular, these satisfy $h^{1,1}(X_i) = h^{2,1}(\tilde{X}_i^{\ast})$ and $h^{2,1}(X_i) = h^{1,1}(\tilde{X}_i^{\ast})$ for $i=0,1,2,3$. 
We refer to \cite{GP, Kan, Kap_prim, Kap_tran} for the Hodge numbers of $X_0$. 
\end{Prop}

\begin{Pf}
We only calculate the Hodge numbers of $\tilde{X}_1^{\ast}$. 
Since the other cases are similar, we leave them to the readers. 
First, we calculate the Euler number of $\tilde{X}_1^{\ast}$. 
By the properties of the Euler number, we have 
\begin{align*}
e(T_{\bar{0}} \times_{\PP^1} T_1) &= e(\pi^{-1}(0))e(\pi'^{-1}(0)) + e(\pi^{-1}(\infty))e(\pi'^{-1}(\infty)) = 45
\end{align*}
where $\pi : T_{\bar{0}} \rar \PP^1$ and $\pi' : T_1 \rar \PP^1$ are projections. 
Since the resolution $\tilde{X}_1^{\ast} \rar T_{\bar{0}} \times_{\PP^1} T_1$ is small, an ordinary double point of $T_{\bar{0}} \times_{\PP^1} T_1$ is replaced by $\PP^1$. 
Then we have 
\begin{align*}
e(\tilde{X}_1^{\ast}) = e(T_{\bar{0}} \times_{\PP^1} T_1) + 45 = 90. 
\end{align*}
Next, we calculate $h^{1,1}(\tilde{X}_1^{\ast})$. 
We use the following result from \cite[Proposition 7.1]{Scho}. 
\begin{Lem}
Let $S_1$ and $S_2$ be relatively minimal rational elliptic surfaces with sections. 
Let $\pi_i : S_i \rar \PP^1$ be the projection. 
Let $D_i$ be the image of singular fibers under the projection $\pi_i$. 
We define $D' =D_1 \cap D_2$. 
Assume that singular fibers of $S_i$ over $D'$ are of type $\I_b$ for $i=1,2$. 
For $s \in D_i$, we denote by $b_i(s)$ the number of irreducible components of the singular fiber $\pi_i^{-1}(s)$. 
Let $\eta$ be the generic point of $\PP^1$. 
We define $d=1$ if $(S_1)_{\eta}$ and $(S_2)_{\eta}$ are isogenous and $d=0$ if otherwise. 
Then $h^{1,1}$ of Schoen's Calabi--Yau $3$-fold is given by 
\begin{align*}
d + 19 + \sum_{s \in D'} b_1(s)b_2(s) - \sum_{s \in D'} b_1(s) - \sum_{s \in D'} b_2(s) + \# D'. 
\end{align*}
\end{Lem}
\noindent
Applying the above lemma to our case, we obtain $h^{1,1}(\tilde{X}_1^{\ast}) = 47$. 
By the relation $e(\tilde{X}_1^{\ast}) = 2(h^{1,1} - h^{2,1})$, we obtain $h^{2,1}(\tilde{X}_1^{\ast}) = 2$ as desired. 
\qed
\end{Pf}

\section{Families of Schoen's Calabi--Yau $3$-folds}\label{sec5}
In this section, we will construct families of Calabi--Yau $3$-folds of $\tilde{X}_0^{\ast}$ and $\tilde{X}_1^{\ast}$ which are expected to be mirror families of $X_0$ and $X_1$, respectively. 
When constructing these families, we will use the relation between Hadamard products and fiber products in \cite[Section 3]{SvS}. 

\subsection{Mirror families of linear sections of $\Join(\Gr, \Gr)$}\label{mirror_GG}
Let $S(5)$ be the Shioda's modular surface associated with the congruence subgroup $\Gamma_1(5)$. 
Let $\pi:S(5) \rar \PP^1$ be the projection to the modular curve. 
We fix a homogeneous coordinate of $\PP^1$ as done in \cite{Zag}. 
Then the map $\pi$ has $\I_5$-fibers over the points $0$ and $\infty$ and $\I_1$-fibers over the points $[1,z]$ where $z^2 + 11z - 1 = 0$. 
\vspace{1mm}

We set $S_i = S(5)$ and $S_i \rar \PP^1 \; (i=1,2)$ as above. 
Let $\mu : \PP^1 \times \PP^1 \dashrightarrow \PP^1$ be the multiplication map defined by 
\begin{align*}
([s_0, s_1], [t_0, t_1]) \mapsto [s_0 t_0, s_1 t_1]. 
\end{align*}
The map $\mu$ is defined outside of $B = \{(0, \infty), (\infty, 0)\}$. 
Let $\widetilde{\PP^1 \times \PP^1}$ be the blow-up of $\PP^1 \times \PP^1$ along $B$. 
Then the rational map $\mu$ extends to the morphism $\widetilde{\PP^1 \times \PP^1} \rar \PP^1$. 
We define $\widetilde{S_1 \times S_2}$ by the fiber product of $\widetilde{\PP^1 \times \PP^1} \rar \PP^1 \times \PP^1$ and $S_1 \times S_2 \rar \PP^1 \times \PP^1$: 
\[ 
\begin{tikzcd}
\widetilde{S_1 \times S_2} \arrow{d} \arrow{r} & S_1 \times S_2 \arrow[d] \\ 
\widetilde{\PP^1 \times \PP^1} \arrow{r} & \PP^1 \times \PP^1
\end{tikzcd} 
\]
We define a morphism $p : \widetilde{S_1 \times S_2} \rar \PP^1$ by the composition of the projection $\widetilde{S_1 \times S_2} \rar \widetilde{\PP^1 \times \PP^1}$ and the morphism $\widetilde{\PP^1 \times \PP^1} \rar \PP^1$. 

\begin{Prop}
Let $D_0 = \{ 0, \infty \}$ and $D_1 = \{[1, z] \mid (z+1)(z^2 - 123 z + 1) = 0\}$. 
For $[z_0, z_1] \in \PP^1 \setminus (D_0 \cup D_1)$, the fiber of the map $p : \widetilde{S_1 \times S_2} \rar \PP^1$ over $[z_0, z_1]$ is isomorphic to $X_0^{\ast}$ constructed in Section \ref{CY_from_modular}. 

\end{Prop}

\begin{Pf}
For each $[1,z] \in \PP^1 \setminus D_0$, the fiber $p^{-1}([1,z])$ is the fiber product of rational elliptic surfaces which are isomorphic to $S(5)$. 
Indeed, the fiber  $p^{-1}([1,z])$ is isomorphic to the fiber product $S(5) \times_{\PP^1} S(5)'$ where $S(5)'$ is isomorphic to $S(5)$ but a morphism $S(5)' \rar \PP^1$ is given by the composition of $\pi:S(5) \rar \PP^1$ and the map $\PP^1 \rar \PP^1$ where  
\begin{align*}
[s_0, s_1] \mapsto [s_1,  z s_0]. 
\end{align*}
The condition that singular fibers of $S(5) \rar \PP^1$ and $S(5)' \rar \PP^1$, other than over $0$ and $\infty$, are located over the different positions is equivalent to $[1, z] \notin D_1$. 
\qed
\end{Pf}

For $[z_0, z_1] \in \PP^1 \setminus (D_0 \cup D_1)$, the fiber $p^{-1}([z_0, z_1])$ has a small resolution, which is a Calabi--Yau $3$-fold $\tilde{X}^{\ast}_0$. 
Therefore we obtain a family of Schoen's Calabi--Yau $3$-folds over $\PP^1$. 
For $[z_0, z_1] \in D_1$, the fiber $p^{-1}([z_0, z_1])$ has an extra nodal point, so $D_1$ is a discriminant locus of the family. 
\vspace{3mm}

This family of Calabi--Yau $3$-folds have two large complex structure limit points at $0$ and $\infty$. 
This result reflects the fact that there exists Calabi--Yau $3$-folds $X$ and $Y$ which are derived equivalent \cite{KP} but not isomorphic to each other \cite{BCP, OR}. 
Note that both Calabi--Yau $3$-folds are linear sections of $\Join(\Gr, \Gr)$. 

\begin{Rem}
A mirror family of linear section Calabi--Yau $3$-folds in $\Join(\Gr, \Gr)$ was already studied by \cite{Kap, Miu}. 
Their constructions are based on the toric degeneration of $\Join(\Gr, \Gr)$ and conifold transitions of linear section Calabi--Yau $3$-folds. 
\end{Rem}

\subsection{Mirror families of linear sections of $\Join(\Gr, \PP^2 \times \PP^2)$}\label{mirror_G(2,V_5)_P2P2}
First we construct a family of elliptic curves whose period integrals coincide with the $I$-function of linear sections of $\PP^2 \times \PP^2$ (see Appendix \ref{I_X_1}). 
We denote by $\PP^1_{\bm{s}}, \PP^2_{\bm{t}}, \cdots$, etc, the projective spaces with homogeneous coordinates $\bm{s}=[s_0, s_1], \bm{t}=[t_0, t_1, t_2]$, and so on. 

\begin{Prop}\label{mirrorP2P2}
There exists a three dimensional projective variety ${\tt{T}}$ and a morphism $\pi : {\tt{T}} \rar \PP^2_{\bm{t}}$ such that 
\begin{enumerate}
\item A general fiber of $\pi$ is an elliptic curve. 
\item For $[z_1, z_2] \in \PP^1$, let $\ell_{[z_1,z_2]}$ be a line in $\PP^2_{\bm{t}}$ defined by $\{[t_0, t_1, t_2] \in \PP^2_{\bm{t}} \mid z_2 t_1 - z_1 t_2 = 0\}$. Then the inverse image $\pi^{-1}(\ell_{[z_1, z_2]})$ is either the Shioda's modular surface $S(6)$ associated with the congruence subgroup $\Gamma_1(6)$ if $[z_1, z_2] = [1,1]$ or a rational elliptic surface with section if $[z_1, z_2] \neq [1,0], [1,-1], [1,1],[0,1]$. 
The latter case has five singular fibers. 
\item Let us fix $[z_1, z_2] \neq [1,0], [1,-1], [0,1]$. 
For the map $\pi^{-1}(\ell_{[z_1, z_2]}) \rar \ell_{[z_1,z_2]}$, the fiber over the point $[1,0,0]$ is of type $\I_6$ and the fiber over the point $[0, z_1, z_2]$ is of type $\I_3$. 
\end{enumerate}
\end{Prop}

\begin{Pf}
We will sketch our construction of ${\tt{T}}$. %
Let us consider the following complete intersection ${\tt{T}}_0$ in $\PP^2_{\bm{x}} \times \PP^2_{\bm{y}} \times \PP^2_{\bm{c}}$: 
\begin{align}\label{sp_family}
&c_0 x_0 y_0 + c_1 x_1 y_2 + c_2 x_2 y_1 = 0, \hspace{3mm} c_0 x_1 y_1 + c_1 x_2 y_0 + c_2 x_0 y_2 = 0, \\
&c_0 x_2 y_2 + c_1 x_0 y_1 + c_2 x_1 y_0 = 0. \notag
\end{align}
There is a natural projection $\pi_0 : {\tt{T}}_0 \rar \PP^2_{\bm{c}}$ coming from the projection to the third factor. 
The variety ${\tt{T}}_0$ is smooth outside of $81$ ordinary double points, which are on the fibers of $\pi_0$ over $\{[1, -\zeta^i, 0], [1, 0, -\zeta^j], [0,1,-\zeta^k] \mid 0 \leq i,j,k \leq 2 \}$ where $\zeta = e^{2 \pi i / 3}$. 

Let $G \cong \ZZ_3^{\pl 3}$ be the subgroup of $\mathrm{Aut}(\PP^2_{\bm{x}} \times \PP^2_{\bm{y}} \times \PP^2_{\bm{c}})$ generated by 
\begin{align*}
([x_0,x_1,x_2],[y_0,y_1,y_2],[c_0,c_1,c_2]) &\mapsto ([x_0,\zeta x_1,x_2],[y_0,y_1, \zeta y_2],[c_0,\zeta c_1,c_2]), \\
([x_0,x_1,x_2],[y_0,y_1,y_2],[c_0,c_1,c_2]) &\mapsto ([x_0,x_1,\zeta x_2],[y_0,\zeta y_1,y_2],[c_0,c_1,\zeta c_2]), \\
([x_0,x_1,x_2],[y_0,y_1,y_2],[c_0,c_1,c_2]) &\mapsto ([x_0,\zeta x_1,\zeta^2 x_2],[y_0,\zeta y_1, \zeta^2 y_2],[c_0,c_1,c_2]). 
\end{align*}
The group $G$ acts on ${\tt T}_0$. 
We define ${\tt T}$ by the quotient ${\tt T}_0/G$ of ${\tt T}_0$ by $G$. 
Then the projection $\pi_0 : {\tt T}_0 \rar \PP^2_{\bm{c}}$ naturally induces the morphism $\pi : {\tt T} \rar \PP^2_{\bm{c}}/G =: \PP^2_{\bm{t}}$ where $[t_0, t_1, t_2] = [c_0^3, c_1^3, c_2^3]$. %

It is straightforward to show that $\pi^{-1}(\ell_{[1,1]})$ is isomorphic to the modular surface $S(6)$ by comparing the functional invariants and also homological invariants of $\pi^{-1}(\ell_{[1,1]})$ with those of $S(6)$. 
\qed
\end{Pf}

\begin{Rem}
The above special family of equations (\ref{sp_family}) has been chosen so that the period integrals of the family of elliptic curves coincide with the period integrals of toric mirror construction of the complete intersection $(1,1) \cap (1,1) \cap (1,1)$ in $\PP^2 \times \PP^2$ by Batyrev--Borisov \cite{BB} (also \cite{HT_det}). 
\end{Rem}

We will construct a family of Schoen's Calabi--Yau $3$-folds $\tilde{X}^{\ast}_1$. %
Consider the above family ${\tt T} \rar \PP^2_{\bm{t}}$ and also $S := S(5) \rar \PP^1_{\bm{s}}$. %
By the analogy of the construction in Section \ref{mirror_GG}, we first consider a rational map $\mu : \PP^1_{\bm{s}} \times \PP^2_{\bm{t}} \dashrightarrow \PP^2_{\bm{z}}$ given by 
\begin{align*}
([s_0, s_1] , [t_0, t_1, t_2]) \mapsto [s_0 t_0, s_1 t_1, s_1 t_2]. 
\end{align*}
This map $\mu$ is defined outside of $B_1 \sqcup B_2$ where 
\begin{align*}
B_1 &= \{([1,0], [0, t_1, t_2]) \mid (t_1, t_2) \in \CC^2 \setminus \{(0,0)\} \}, \\
B_2 &= \{([0,1], [1, 0, 0]) \}. 
\end{align*}
Let $\widetilde{\PP^1_{\bm{s}} \times \PP^2_{\bm{t}}}$ be the blow-up of $\PP^1_{\bm{s}} \times \PP^2_{\bm{t}}$ along $B_1 \sqcup B_2$. 
Then the rational map $\mu$ extends to the morphism $\widetilde{\PP^1_{\bm{s}} \times \PP^2_{\bm{t}}} \rar \PP^2_{\bm{z}}$. 
We define $\widetilde{S \times {\tt T}}$ by the fiber product expressed in the following diagram: 
\[ 
\begin{tikzcd}
\widetilde{S \times {\tt T}} \arrow{d} \arrow{r} & S \times {\tt T} \arrow[d] \\ 
\widetilde{\PP^1_{\bm{s}} \times \PP^2_{\bm{t}}} \arrow{r} & \PP^1_{\bm{s}} \times \PP^2_{\bm{t}}
\end{tikzcd} 
\]
We define a morphism $p : \widetilde{S \times {\tt T}} \rar \PP^2_{\bm{z}}$ by the composition of the projection $\widetilde{S \times {\tt T}} \rar \widetilde{\PP^1_{\bm{s}} \times \PP^2_{\bm{t}}}$ and the morphism $\mu : \widetilde{\PP^1_{\bm{s}} \times \PP^2_{\bm{t}}} \rar \PP^2_{\bm{z}}$. %
Note that $\mu^{-1}([z_0, z_1, z_2]) \cong \PP^1$ for general $[z_0, z_1, z_2] \in \PP^2_{\bm{z}}$. %

\begin{Prop}
Let $D_2 = \{dis(z_0,z_1,z_2) = 0\}$ be a divisor in $\PP^2_{\bm{z}}$ where 
\begin{align*}
dis(z_0,z_1,z_2) = \prod_{\zeta^2+11\zeta-1=0} ((\zeta z_0 + z_1 + z_2)^3 - 27 \zeta z_0 z_1 z_2). 
\end{align*}
We also define $D_0 = \{z_0 z_1 z_2 = 0\}$ and $D_1 = \{z_1 + z_2 = 0 \}$. 
Then for any $[z_0, z_1, z_2] \in \PP^2_{\bm{z}} \setminus (D_0 \cup D_1 \cup D_2)$, a fiber of the morphism $p : \widetilde{S \times {\tt T}} \rar \PP^2_{\bm{z}}$ over $[z_0, z_1, z_2]$ is the fiber product of rational elliptic surfaces with sections. %
A crepant resolution of the fiber $p^{-1}([z_0, z_1, z_2])$ is deformation equivalent to Schoen's Calabi--Yau $3$-folds $\tilde{X}^{\ast}_{1}$ constructed in Section \ref{CY_from_modular}. 
\end{Prop}

\begin{Pf}
Let us take a point $[z_0, z_1, z_2] \in \PP^2_{\bm{z}} \setminus (D_0 \cup D_1)$. 
A fiber $p^{-1}([z_0, z_1, z_2])$ is isomorphic to $S(5) \times_{\PP^1_{\bm{s}}} {\tt T}_{z_1, z_2}$ where ${\tt T}_{z_1,z_2}$ is the inverse image of $\{z_2 t_1 - z_1 t_2 = 0\} \subset \PP^2_{\bm{t}}$ under the map $\pi : {\tt T} \rar \PP^2_{\bm{t}}$. 
Note that we identify $\PP^1_{\bm{s}} \cong \{ z_2 t_1 - z_1 t_2 = 0 \} \subset \PP^2_{\bm{t}}$ by the embedding: %
\begin{align*}
[s_0, s_1] \mapsto [z_0 s_1, z_1 s_0, z_2 s_0]. 
\end{align*}
By Proposition \ref{mirrorP2P2}, if $z_1 = z_2$ then ${\tt T}_{z_1, z_2}$ is a modular surface $S(6)$, and if $z_1 \neq \pm z_2$ then ${\tt T}_{z_1,z_2}$ is an elliptic surface with five singular fibers. 
We define by $F_1 \subset \PP^1_{\bm{s}}$ and $F_2 \subset \PP^1_{\bm{s}}$ the images of the singular fibers of the elliptic fibrations $S(5) \rar \PP^1_{\bm{s}}$ and ${\tt T}_{z_1, z_2} \rar \PP^1_{\bm{s}}$, respectively. 
Then the condition $F_1 \cap F_2 = \{[1,0], [0,1]\}$ is equivalent to the condition $[z_0, z_1, z_2] \notin D_2$. 
Hence for $[z_0, z_1, z_2] \in \PP^2_{\bm{z}} \setminus (D_0 \cup D_1 \cup D_2)$ satisfying $z_1 = z_2$, the fiber $p^{-1}([z_0, z_1, z_2])$ describes $X_1^{\ast}$ in Conjecture \ref{main_Conj}. 
Therefore the claim follows. 
\qed 
\end{Pf}

Let $[1, -z_1, -z_2]$ be the affine coordinate of $\PP^2_{\bm{z}}$ around the point $[1,0, 0]$. 
By the construction of the above family of Calabi--Yau $3$-folds, there exists a period $\omega_0(z_1, z_2)$ given by the Hadamard product of period integrals of two families of elliptic curves. 
The explicit form of $\omega_0$ is given in (\ref{Hadam1}) where we set $x_i=z_i$. 
\vspace{3mm}

It is straightforward to compute the Picard--Fuchs operators which annihilate $\omega_0(z_1,z_2)$. 
They are given by 
\begin{align*}
P_1 &= p_{20} \theta_1^2 + p_{11} \theta_1 \theta_2 + p_{02} \theta_2^2 + p_{10} \theta_1 + p_{01} \theta_2 + p_{00}, \\
P_2 &= z_1 \theta_2^3 - z_2 \theta_1^3
\end{align*}
where we set $\theta_i = z_i \frac{\partial}{\partial z_i}$. 
The coefficients of $\theta_1^i \theta_2^j \; (0 \leq i + j \leq 2)$ in $P_1$ are 
\begin{align*}
p_{20} &= z_1^2 + 11z_1 z_2 + 10 z_2^2 + 11z_1 + 11z_2 -1, \\
p_{11} &= 5z_1^2 + 10 z_1 z_2 + 5z_2^2 + 22z_1 + 22z_2 + 1, \\
p_{02} &= 10z_1^2 + 11z_1 z_2 + z_2^2 + 11z_1 + 11z_2 -1, \\
p_{10} &= (z_1 + z_2) (2z_1 + 5z_2 + 11), \\
p_{01} &= (z_1 + z_2) (5z_1 + 2z_2 + 11), \\
p_{00} &= (z_1 + z_2) (z_1 + z_2 + 3). 
\end{align*}
Let us study these differential operators over $\PP^2_{\bm{z}}$. 
Around the points $[0,1,0]$ and $[0,0,1]$, we have the following facts. 

\begin{Prop}
Let $[w_0, 1, -w_2]$ be the affine coordinate of $\PP^2_{\bm{z}}$ around the point $[0,1,0]$. 
Let $Q_1$ and $Q_2$ be the following gauge transforms of the operators $P_1, P_2$ around the point $[0,1,0]$: %
\begin{align*}
Q_1 = w_0 P_1(w_0, w_2) w_0^{-1}, Q_2 = P_2(w_0, w_2). 
\end{align*}
Then $Q_1$ and $Q_2$ determine the period integrals of the above family around the point $[0,1,0]$. 
Similar results hold for the point $[0,0,1]$. 
\end{Prop}
\begin{Pf}
This is a standard fact. We refer \cite[Proposition 6.1]{HT_det}, for example, for details. %
\qed
\end{Pf}

\begin{Prop}
There exists a period integral which is holomorphic around the point $[0, 1, 0]$. 
More precisely, it is given by 
\begin{align*}
\sum_{d_0', d_2' \geq 0} \left( \sum_{k=0}^{d_0'}\binom{d_0'}{k}^2\binom{d_0' + k}{k} \right) \frac{(d_0' + d_2')!^3}{(d_{0}'!)^3 (d_{2}'!)^3} w_0^{d_0'} w_2^{d_2'}
\end{align*}
where $[w_0, 1, -w_2]$ is the affine coordinate of $\PP^2_{\bm{z}}$. 
Similar results hold for the point $[0,0,1]$. 
\end{Prop}
\begin{Pf}
This follows by the argument in \cite[Section 3]{SvS}. %
\qed
\end{Pf}

By calculating the monodromies around the coordinate hyperplanes of $\PP^2_{\bm{z}}$, we obtain the following result. 
\begin{Prop}
The points $[1,0,0], [0,1,0]$, and $[0,0,1]$ are large complex structure limit points in the sense of Morrison \cite{Mor}. 
\end{Prop}

We can calculate the $g=0$ Gromov--Witten invariants for $X_1$ (resp.\ $Y_1$) by applying J.\ Brown's result \cite{Bro} or the abelian/nonabelian correspondence \cite{BCFK} to resolved joins $\J_1$ and $\cO(L)^{\pl 8}$ (resp.\ $\J'_1$ and $\cO(L')^{\pl 11}$). 
The tables of BPS numbers of $X_1$ and $Y_1$ are shown in Appendix \ref{appA}. 

\begin{Rem}
We can extract enumerative invariants from Picard--Fuchs operators for each large complex structure limit points $[1,0,0], [0,1,0]$, and $[0,0,1]$ (see \cite{CK} and references therein for the details). 
By construction, the invariants obtained from the point $[1,0,0]$ coincide with the $g=0$ Gromov--Witten invariants of $X_1$. 
By computation, we can check that the invariants obtained from the points $[0,1,0]$ and $[0,0,1]$ coincide with the $g=0$ Gromov--Witten invariants of $Y_1$ for sufficiently large degree. 
\end{Rem}

\appendix
\section{Lists of enumerative invariants}\label{appA}
In this Appendix, we give the lists of the BPS numbers obtained from the $g=0$ Gromov--Witten invariants of $X_1$ and $Y_1$. 
\vspace{3mm}

\subsection{} 
Let $X_1$ be the Calabi--Yau $3$-fold constructed in Section \ref{The_case_del_Pezzo}. 
We fix a basis of $H^2(X_1, \QQ)$ as $D_1, D_2$ where $D_i$ is the restriction of the pull-back of the hyperplane class of the $i$-th factor of $\PP^2 \times \PP^2$. 
We represent a homology class $C \in H_2(X_1, \ZZ)$ by the intersection numbers $(d_1, d_2):=(D_1.C, D_2.C)$. 

\begin{table}[H]\centering\footnotesize
 \begin{tabular}{c |l l l l l l}
\hline
 $d_1 \setminus d_2$ & $0$ & $1$ &$2$ &$3$&$4$&$5$ \\
\hline
0& 0 & 120 & 105 & 105 & 120 & 90 \\
1& 120 & 2085 & 15690 & 83400 & 362850 & 1365060 \\
2& 105 & 15690 & 569475 & 9690270 & 107459880 & 901887570 \\
3& 105 & 83400 & 9690270 & 418812780 & 10086474180 & 164859436335 \\
4& 120 & 362850 & 107459880 & 10086474180 & 472152998265 & 13800385325580 \\
5& 90 & 1365060 & 901887570 & 164859436335 & 13800385325580 & 675995017391805 \\
6& 120 & 4621020 & 6204484125 & 2041590595410 & 286700834960805 & 22351196770131870 \\
7& 105 & 14399490 & 36701125005 & 20496053409240 & 4593254607725475 & 546563929916334210 \\
8& 105 & 41932200 & 192593575110 & 174405931797135 & 59937858896889555 & 10518492857890739820 \\
9& 120 & 115485075 & 916315955820 & 1297448843314125 & 661998422042833065 & 166511015537610566130 \\
10& 90 & 303166710 & 4015843886955 & 8630138044756890 & 6364684023911207415 & 2240097475662256021890
\end{tabular}
\caption{BPS numbers $n_0^{X_1}(d_1,d_2)$ of $X_1$}
\end{table}

\subsection{}
Let $Y_1$ be the Calabi--Yau $3$-fold constructed in Section \ref{The_case_del_Pezzo}. 
We fix a basis of $H^2(Y_1, \QQ)$ as $L', D'$ where $L'$ is the restriction of the relative hyperplane class of the resolved join and $D'$ is the pull-back of the hyperplane class of $\PP^2$ by the map of the elliptic fibration on $Y_1$. 
We represent a homology class $C' \in H_2(Y_1, \ZZ)$ by the intersection numbers $(d'_1, d'_2):=(L'.C', D'.C')$. 

\begin{table}[H]\centering\footnotesize
 \begin{tabular}{c |l l l l l l}
\hline
 $d'_1 \setminus d'_2$ & $0$ & $1$ &$2$ &$3$&$4$&$5$ \\
\hline
0&0&30&0&0&0&0\\
1&105&330&105&0&0&0\\
2&120&2865&6585&2865&120&0\\
3&120&17400&151260&283755&151260&17400\\
4&105&87150&2141265&11044335&18347055&11044335\\
5&90&368670&22279830&256967580&974066175&1488072900\\
6&105&1377840&186120810&4267143150&31595446320&97322962410\\
7&120&4644030&1311908070&55405726800&729262582320&4007703642030\\
8&120&14441100&8065898475&594374999280&13050194338080&118409369639565\\
9&105&42003450&44272540830&5463083502630&191094069663765&2712537543756540\\
10&90&115593255&220759120890&44140588111590&2375090868607470&50686607599977960
\end{tabular}
\caption{BPS numbers $n_0^{Y_1}(d'_1,d'_2)$ of $Y_1$}
\end{table}

\section{$I$-functions of complete intersections of resolved joins}
We calculate $I$-functions of Calabi--Yau $3$-folds $X_0$ and $X_1$ discussed in Section \ref{mirror_GG} and \ref{mirror_G(2,V_5)_P2P2}. 
We refer to \cite{CG} for definitions and properties of $I$-functions. 
Note that the $I$-functions of Calabi--Yau $3$-folds are related to period integrals of mirror families. 

\subsection{The $I$-function of $X_0$}
To calculate the $I$-function of $X_0$, we use J.\ Brown's formula about $I$-functions of toric bundles. 
Let us recall that $X_0$ is a complete intersection of ten relative hyperplanes in the resolved join $\PP_{\Gr \times \Gr}(\cO(-H_1) \pl \cO(-H_2))$. 

The following proposition is an immediate consequence of J.\ Brown's formula. 
This result is already mentioned by S.\ Galkin in his talk \cite{Gal}. 

\begin{Prop}\label{I-function_X_0}
Let $J^{\Gr}(\tau,z) = e^{\tau H/z} \sum_{d=0}^{\infty} J^{\Gr}_{d}(\tau,z) Q^d$ be the $J$-function of $\Gr$ where $H$ is the Schubert divisor of $\Gr$. 
We denote by $G(2, V_5)_i$ the $i$-th factor of $G(2, V_5) \times G(2, V_5)$. 
Then the $I$-function of $X_0$ is 
\begin{align*}
I_{X_0}(\tau, z) = e^{\tau L / z} \sum_{d = 0}^{\infty} \prod_{1 \leq i \leq 2}\prod_{m=1}^{d}(L + m z)^5 J^{\Gr_i}_d(\tau/2, z)|_{X_0} Q^d
\end{align*}
where $J^{\Gr_i}_d(\tau/2, z)|_{X_0}$ is the restriction of cohomology classes of $\Gr_i$ to $X_0$. 
\end{Prop}

\begin{Pf}
Applying J.\ Brown's formula \cite{Bro} and quantum Lefschetz theorem \cite{CG} to a complete intersection $X_0$ of the resolved join $\PP := \PP_{\Gr \times \Gr}(\cO(-H_1) \pl \cO(-H_2))$, we have the twisted $I$-function 
\begin{align*}
I({\bm \tau}, z) = e^{{\bm \tau} /z} \sum_{d_1, d_2, d_3 \geq 0} \frac{\prod_{m=1}^{d_3} (L + mz)^{10}}{\prod_{1 \leq i \leq 2}\prod_{m=1}^{d_3-d_i}(L - H_i + mz)}J^{G(2, V_5)_1}_{d_1}(\tau_1, z) J^{G(2, V_5)_2}_{d_2}(\tau_2, z) e^{d_3 \tau_3} Q_1^{d_1} Q_2^{d_2} Q_3^{d_3}
\end{align*}
where ${\bm \tau} = \tau_1 H_1 + \tau_2 H_2 + \tau_3 L$. 
Since $(L - H_i)|_{X_0} \sim 0$, we have an injective homomorphism of $\mathrm{NE}(X_0)$ to $\mathrm{NE}(\PP)$ and the injective homomorphism of Novikov rings given by $Q \mapsto Q_1 Q_2 Q_3$ where $Q$ is a generator of the Novikov ring of $X_0$. 
Passing from the twisted $I$-function to the $I$-function of $X_0$, we obtain the desired form. 
\qed
\end{Pf}

\begin{Rem}
The coefficient of $z^0$ in the $I$-function of $X_0$ is 
\begin{align*}
\sum_{d=0}^{\infty} \left( \sum_{k=0}^{d}\binom{d}{k}^2 \binom{d+k}{k} \right)^2 x^d 
\end{align*}
where $x = e^{\tau} Q$. 
This is a Hadamard product of the period integrals of Shioda's modular surface associated with $\Gamma_1(5)$ (see, for example, \cite{Zag}). 
\end{Rem}

\subsection{The $I$-function of $X_1$}\label{I_X_1}
Similar to $X_0$, we can calculate the $I$-function of $X_1$. 
Let us recall that $X_1$ is a complete intersection of eight relative hyperplanes in the resolved join $\PP_{\Gr \times \PP^2 \times \PP^2}(\cO(-H_1) \pl \cO(-H_2))$ where $H_1$ is the pull-back of the Schubert divisor class of $\Gr$ and $H_2$ is the pull-back of the $(1,1)$ class of $\PP^2 \times \PP^2$. 

\begin{Prop}
Let $J^{\Gr}(\tau, z)$ be the $J$-function of $\Gr$ as above. 
Let $J^{\PP^2 \times \PP^2}(\tau_1, \tau_2, z) = e^{\tau_1 D_1+ \tau_2 D_2/z} \sum_{d_1, d_2 \geq 0} J^{\PP^2 \times \PP^2}_{d_1, d_2}(\tau_1, \tau_2, z) Q_{1}^{d_1} Q_{2}^{d_2}$ be the $J$-function of $\PP^2 \times \PP^2$ where $D_i \; (i=1,2)$ is the hyperplane class of the $i$-th factor of $\PP^2 \times \PP^2$. 
Then the $I$-function of $X_1$ denoted by $I_{X_1}(\tau_1, \tau_2, z)$ is 
\begin{align*}
e^{\tau_1 D_1 + \tau_2 D_2/z} \sum_{d_1, d_2 \geq 0} 
\prod_{m=1}^{d_1+d_2} (L + m z)^5 J^{\Gr}_{d_1 + d_2}(0, z)|_{X_1}
\prod_{m=1}^{d_1 + d_2} (L + m z)^3 J^{\PP^2 \times \PP^2}_{d_1, d_2}(\tau_1, \tau_2, z)|_{X_1} Q_1^{d_1} Q_2^{d_2}
\end{align*}
where $D_i$ is the divisor class on $X_1$ corresponding to the elliptic fibrations on $X_1$. 
\end{Prop}

\begin{Pf}
The proof is similar to Proposition \ref{I-function_X_0}. We omit the details. 
\qed
\end{Pf}

\begin{Rem}
The coefficient of $z^0$ in the $I$-function of $X_1$ is 
\begin{align}\label{Hadam1}
\sum_{d_1, d_2 \geq 0} \left( \sum_{k=0}^{d_1 + d_2}\binom{d_1 + d_2}{k}^2 \binom{d_1 + d_2 + k}{k} \right) \frac{(d_1 + d_2)!^3}{(d_{1}!)^3 (d_{2}!)^3} x_1^{d_1} x_2^{d_2}
\end{align}
where $x_1 = e^{\tau_1} Q_1$ and $x_2 = e^{\tau_2} Q_2$. 
This is a generalization of the usual Hadamard product to the case between the power series of one variable and the power series of two variables. 
\end{Rem}

\end{document}